\documentclass[11pt]{article}
\usepackage{amssymb}
\usepackage{amsmath}
\usepackage{amsthm}

\usepackage[margin=1 in]{geometry}

\newcommand{\proofoflem}{\paragraph{\it Proof of lemma.}}

\theoremstyle{plain}
\newtheorem{thrm}{Theorem}
\newtheorem{lem}{Lemma}
\newtheorem{cor}{Corollary}
\newtheorem{prop}{Proposition}
\newtheorem{claim}{Claim}

\theoremstyle{definition}
\newtheorem{defn}{Definition}
\newtheorem{exam}{Example}

\newtheorem*{remark}{Remark}

\makeatletter
\renewcommand\section{\@startsection {section}{1}{\z@}%
                                   {-3.5ex \@plus -1ex \@minus -.2ex}%
                                   {2.3ex \@plus.2ex}%
                                   {\normalfont\large\bfseries}}
\renewcommand\subsection{\@startsection{subsection}{2}{\z@}%
                                     {-3.25ex\@plus -1ex \@minus -.2ex}%
                                     {1.5ex \@plus .2ex}%
                                     {\normalfont\normalsize\bfseries}}
\renewcommand\subsubsection{\@startsection{subsubsection}{3}{\z@}%
                                     {-3.25ex\@plus -1ex \@minus -.2ex}%
                                     {1.5ex \@plus .2ex}%
                                     {\normalfont\normalsize\bfseries\slshape}}
\makeatother

\title{Knot concordance, the point class in instanton homology and  Donaldson invariants}
\author{Yuhan Lim}
\date{Version: February 2024}
\begin{document}

\maketitle

\noindent[\textit{We define an invariant ${\varphi}$ for knots in the $3$-sphere by means of Donaldson invariants and Floer's instanton homology. Some basic properties of this invariant are established and it is shown that ${\varphi}$ coincides with a special case of an invariant defined by Fr{\o}yshov}.]

\section{Introduction}

Let $K$ be a knot in $S^3$. A standard construction to associate a 4-manifold to $K$ is to attach a 2-handle to the boundary of the 4-ball along $K$. This construction is well defined up to diffeomorphism once we specify the attaching map by a given framing of $K$. Let us denote the resultant 4-manifold by $Z_K$, when we use the $0$-framing on $K$. The aim of this paper is to examine some properties of the Donaldson invariants for $Z_K$. These will in turn determine invariants for the knot $K$. 

Since $Z_K$ is a manifold with boundary, the Donaldson invariant will take values in Floer's instanton homology of the boundary 3-manifold. The boundary is in this case $K_0$, the 3-manifold obtained by $0$-surgery on $K$. 
Specifically, over $Z_K$ there is a unique principle $SO(3)$-bundle $P'$ with $w_2$-class $w$ such that $w$ restricts non-trivially over $K_0$. The zero dimensional Donaldson invariant associates a vector
\[
\psi_K\in I^w_*(K_0),
\] 
where the instanton homology groups $I^w_*(K_0)$ are affinely $\mathbb{Z}/8$-graded. We shall work with  rational coefficients as we will be only interested in these groups as vector spaces. 
A simple observation is that $\psi_K=0$ if $K$ is a smoothly slice knot. 
This is because if $K$ is slice then $Z_K$ contains an embedded 2-sphere $S$ with $S\cdot S=0$ and  $w_2(S)\neq 0$. It is well known that this forces the associated Donaldson invariant to vanish. Thus  $\psi_K$ is an obstruction to (smooth) sliceness for $K$. 

In fact $\psi_K$ is an obstruction to $K$ bounding an immersed disk in $B^4$ with only transverse negative self-intersections. This can be seen by resolving the negative self-intersections by taking the connect sum of $Z_K$ with copies of $\overline{\mathbb{CP}^2}$. Choosing $w_2(P')$ to be trivial on the $\overline{\mathbb{CP}^2}$ summands makes no change to the zero-dimensional Donaldson invariant $\psi_K$ (for example see \cite{FM}, 6.6). As before we  obtain an embedded disk $S$ with $S\cdot S=0$ and $w_2(S)\neq 0$ in $Z_K\# n\overline{\mathbb{CP}^2}$,  for $n \gg 1$.

To obtain a finer invariant we will incorporate the action of the $u$-map on $I^w_*(K_0)$. This action is defined to be  the $\mu$-map of the point class. The $u$-map acts by degree $-4$ on $I^w_*(K_0)$. 
Let us set 
\[
{\varphi}(K)= \text{dim span}\{ \psi_K,u^2\psi_K, u^4\psi_K,\dots, u^{2i}\psi_K,\dots\}.
\]

In \cite{Fy} it is shown that for 3-manifolds $Y$ with $b_1\neq 0$ and a choice of $w_2=w$ that is non-trivial on some embedded oriented surface, $(u^2-4)^n=0$ in $I^w_*(Y)$ for some $n\ge 1$. Therefore ${\varphi}(K)$ can be seen as a lower bound on the order of $(u^2-4)$ for $I^w_*(K_0)$. 
We may alternatively define ${\varphi}(K)$ by means of a filtration. The only possible eigenvalue of $u^2$ is $+4$. Let us consider the increasing finite filtration
\[
\{0\}\subset \text{ker}(u^2-4)\subset\text{ker}(u^2-4)^2\subset\text{ker}(u^2-4)^3\subset\dots
\]
defined on $I^w_*(K_0)$. Then ${\varphi}(K)$ is the smallest value $k\ge 0$ for which $(u^2-4)^k\psi_K=0$.

Let $K'$ and $K$ be oriented knots. Define 
 $K'\ge K$ to mean that there is an oriented immersed cylinder $C=[0,1]\times S^1$ with only transverse negative self-intersection points in $[0,1]\times S^3$ such that $\partial C$ is 
 $\{1\}\times K\cup\{0\}\times -K'$. Our convention is 
$\partial ([0,1]\times Y) = \{1\}\times Y \cup \{0\}\times \overline{Y}$. Roughly speaking $K'\ge  K$ is the notion that $K$ is obtained from $K'$ by changing positive crossings in $K'$ to negative crossings. 
We define $K'\le K$ to mean $K^*\ge K'^*$. We remark that the relation $\ge$ was originally  introduced in Cochran-Gompf \cite{CG}, however our definition is 
slightly more restrictive. $K'\ge K$ in our notation implies $K\le K'$ and $K'\ge K$  as defined in \cite{CG} (but not the converse).

\begin{thrm}\label{thrm-one}
Let $K$ be a knot in $S^3$. The following properties hold for ${\varphi}(K)$.
\begin{enumerate}
\item\label{thrm-item1}
If $K'\ge K$ under some choice of orientations, then ${\varphi}(K')\ge {\varphi}(K)$. In particular ${\varphi}$ is a knot concordance invariant.
\item\label{thrm-item2}
If $\Sigma$ is any immersed surface in $B^4$ of genus $g$ having only transverse negative self-intersection points and $\Sigma$ bounds $K$, then
\begin{displaymath}
 {\varphi}(K)\le \lceil g/2\rceil.
 \end{displaymath}
 $\lceil x \rceil$ is the smallest integer greater or equal to $x$.
\item\label{thrm-item3}
Let $K_{+}$ and $K_{-}$ be knots in $S^3$ that in a projection differ only in a positive and negative crossing. Then
\begin{displaymath}
{\varphi}(K_{-})\leq {\varphi}(K_{+})  \leq {\varphi}(K_{-})+1
\end{displaymath}
\item\label{thrm-item4}
Let $K$, $K'$ and $K''$ be  knots each with ${\varphi}>0$.  
Then\footnote{in the earlier version of this paper an incorrect connect sum bound was determined. Corollaries and examples have been updated to reflect this.}
\[
{\varphi}(K\# K'\#K'')\ge {\varphi}(K)+{\varphi}(K')+{\varphi}(K'')-2n,
\]
where $n \ge 1$ is any integer such that $(u^2-4)^n=0$ in $I^w_*(K_0)$, $I^{w'}_*(K'_0)$ and $I^{w''}_*(K''_0)$. 
\end{enumerate}
\end{thrm}

For the right-hand trefoil ${\varphi}\neq 0$, which was observed by Fintushel-Stern \cite{FS}: a neighborhood of a cusp singularity in an elliptic surface is diffeomorphic to $Z_K$ where $K$ is the right-hand trefoil. Since the 0-dimensional Donaldson  invariant is non-trivial for an elliptic K3 surface with choice of $w_2$ that evaluates non-trivially on each generic fiber, the result follows. The genus bound in the Theorem establishes ${\varphi}\le1$ for the right-hand trefoil and therefore ${\varphi}=1$ for this knot.

The $h$-invariant is an invariant of oriented integral homology spheres, defined by Fr{\o}yshov \cite{Fy}.

\begin{thrm}\label{thrm-three}
Let $K$ be a knot in $S^3$. Then
${\varphi}(K) = h(K^*_{-1})$, the $h$-invariant for the integral homology 3-spheres obtained by $(-1)$-surgery on
$K^*$, the mirror of $K$.
\end{thrm}

We remark that  for any knot $K$, $h(K_{-1})$ must  apriori be non-negative since $h$ enjoys a monotonicity property with respect to  negative definite integral homology cobordisms: $K_{-1}$ bounds the  negative definite 4-manifold obtained by attaching a 2-handle to $\partial B^4$ along $K$ with $(-1)$-framing.  

Theorem~\ref{thrm-three}  
directly implies a version of an inequality proven in \cite{Fy2}, restricted to the case of surgery on a knot in $S^3$.

\begin{cor}[\cite{Fy2}]
Let $K$ be a knot in $S^3$. 
Then
$
\lceil g/2\rceil  \ge h(K^*_{-1})
$,
where $g$ is the genus of any immersed surface with only negative transverse self intersection points in the 4-ball that bounds $K$.
\end{cor}

Some of the results in Theorem~\ref{thrm-one} can be deduced via the equivalence but in the context of knots they are more readily established directly from the definition of ${\varphi}(K)$. The proof of Theorem~\ref{thrm-one} does not use the equivalence. However the equivalence does give us some calculations  directly from those in \cite{Fy}.

\begin{exam}
Let $K$ be a positive knot i.e. non-trivial and admits a projection with only positive crossings. Then 
${\varphi}(K) >0$. This follows from a result in \cite[\S 3]{CG} where it is shown that for a non-trivial  positive knot it is possible to change only positive crossings in a sequence of knot projections to transform it to the right-hand trefoil; that is to say any positive knot $\ge$ right-hand trefoil.
\end{exam}

\begin{exam}
Let $T(k)$ be the positive $(p,q)$-torus knots $T(p, q)$ with $p=2$ and $q=2k-1$, $k\ge 2$. It is well-known that $T^*(p,q)_{-1}$ is diffeomorphic to the Brieskorn homology 3-sphere $B(2, 2k-1, 4k-3)$. 
Then ${\varphi}(T(k)) = h(B(2, 2k-1, 4k-3))=\lfloor k/2\rfloor$, the largest integer less than or equal to $k/2$, $k\ge 2$, as shown in \cite{Fy}. \end{exam}

\begin{exam}
If $K$ is a genus 1 or 2 knot then it follows that $u^2-4=0$ (\cite{Fy}). Assume that $K$ is genus 1 or 2 and ${\varphi}(K)>0$; then  ${\varphi}(K)=1$ by Theorem~\ref{thrm-one}(\ref{thrm-item2}). Let $K, K', K''$ be three such knots; then by Theorem~\ref{thrm-one}(\ref{thrm-item4}) we have ${\varphi}(K\# K' \# K'')>0$.
\end{exam}

One may ask if a parallel construction of ${\varphi}(K)$ by attaching a 2-handle to $B^4$ but with a different framing gives different invariants for $K$. In the case of $(-1)$-framing this gives no new information.

\begin{thrm}\label{thrm-two}
Let $K$ be a  knot in $S^3$ and $Z_{K,-1}$ the negative definite 4-manifold obtained by attaching a 2-handle to $B^4$ along $K$ with $(-1)$-framing. Let ${\psi}_{K,-1}\in I_*(K_{-1})$ be the 0-dimensional Donaldson invariant for $Z_{K,-1}$ where we assume $w_2$ is non-zero. Let
\[
{\varphi}_{-1}(K)=\textup{dim span}\{ \psi_{K,-1},u^2\psi_{K,-1}, u^4\psi_{K,-1},\dots, u^{2i}\psi_{K,-1},\dots\}.
\]
Then ${\varphi}_{-1}(K)={\varphi}(K)$. 
\end{thrm} 

We remark that ${\psi}_{K,-1}$ is actually an element of dimension 7 in $I_*(K_{-1})$, As such, the $u$ map is well defined in that dimension, see \cite{Fy}. (Note that \cite{Fy} works with Floer \textit{cohomology}, and the grading with Floer homology is related via the rule $I^*(Y) =I_{5-*}(Y)$.)

This paper is organized as follows. In Section~\ref{prelim} we gather the basic constructions and notations for defining the invariants. In Section~\ref{proofthrm-one} we proof Theorem~\ref{thrm-one} after recalling a result of Munoz \cite{M}. Sections~\ref{conn} and \ref{h-char} lay out some results on stablization by the Poincare homology sphere as a direct application of Fukaya's connected sum result as presented in \cite{Sca}. This shall be one of our main tools. Section~\ref{h-char} may be of independent interest as we present a 4-dimensional characterization for $h(K_{-1})$ leading up to the completion of the proof of Theorem~\ref{thrm-three} in Section~\ref{complete}. The constructions in this paper parallel those in \cite{Fy} and the proof of Theorem~\ref{thrm-three} relies heavily on results proven there.

\section{Preliminaries and definitions}\label{prelim}

For background information about instanton homology we refer to \cite{BD}, \cite{Fy} and \cite{D}.
Henceforth, when we refer to a knot, the assumption will be that this is a knot in $S^3$ unless otherwise stated.  
In the introduction we defined the manifold $Z_K$ obtained by attaching a 2-handle to $B^4$ along $K$ on the boundary. For the goals of this paper it will be more convenient to work with a dual construction. 
Let $K$ be a knot and $\nu K$ a tubular neighbourhood of $K$ that we identify with $S^1\times D^2$ such that $S^1\times\{\mbox{pt}\}$ corresponds to a meridian of the knot and 
$\{\mbox{pt}\}\times \partial D^2$ with a Seifert longitude. We can think of this as the standard framing of $K$ in $S^3$. Remove $\nu K$ and reglue to the complement in a manner that switches the median and Seifert longitude, this effects $0$-surgery on $K$. The core $S^1\times \{0\}$ in $S^1\times D^2$ now represents a framed knot $K'$ in $K_0$. Attaching a 2-handle to $K_0\times[0,1]$ along $K'\times\{1\}$ with the mentioned framing will result in an oriented cobordism 
\begin{equation}\label{cob-one}
W:K_0\to S^3.
\end{equation} 
To keep track of ingoing and outgoing components of a cobordism we shall adopt the notation $W: Y\to Y'$ to mean $\partial W = Y'\cup \overline{Y}$. We regard $Y$ as the ingoing component and $Y'$ the outgoing.
Closing off the $S^3$ component of $W:K_0\to S^3$ with a 4-ball we obtain a cobordism
\begin{displaymath}
W_K: K_0 \to \emptyset.
\end{displaymath}
This is our first construction. Our second construction is to  attach a 2-handle to the 4-ball along $K$ using the standard 0-framing on $K$. This gives us
\begin{displaymath}
W'_K: \emptyset\to K_0.
\end{displaymath}
$W'_K$ is exactly the manifold $Z_K$ mentioned in the introduction.
The two constructions $W_K$ and $W_K'$ are inverses of each other except for an additional reflection. If we reverse the roles of ingoing and outgoing components of $W:Y\to Y'$ we get $W^*:\overline{Y}'\to \overline{Y}$. Crucial to this operation is that the orientation of $W$ is unchanged. We have the relations:
\begin{displaymath}
(W_K)^* = W'_{K^*} , \quad (W'_K)^* = W_{K^*}
\end{displaymath}
where $K^*$ denotes the mirror of $K$.

Geometrically the unique choice of $w$ for $K_0$  can be represented by a meridian of $\nu K$ thinking of it as a curve in the complement of $K$ and hence $K_0$. This class extends over $W_K$: geometrically consider the embedded disk $D$ that is dual to the core of the attached 2-handle. $\partial D$ is a embedded curve in $K_0$ that is a meridian for $\nu K$. An identical construction can be carried out for $W'_K$.  With these choices of $w=w_2$ the cobordisms $W_K$, $W'_K$ define 
\begin{displaymath} 
\delta_K: I^w_*(K_0)\to\mathbb{Z},\quad \delta'_K \in I^w_*(K_0)
\end{displaymath}
where $\delta_K$ is a homomorphism. $\delta'_K$ is what was denoted by $\psi_K$ in the introduction:  we change notation to be consistent with constructions in \cite{Fy}.
To be specific extend $W_K$, $W'_K$ to be cylindrical-end manifolds $\widehat{W}_K$, $\widehat{W}'_K$ by appending semi-infinite tubes $(-\infty,0]\times K_0$, $[0, \infty)\times K_0$ to the ingoing and outgoing boundary components, respectively. Similarly extend any $SO(3)$-bundle $P$ on the compact 4-manifolds with boundary to be compatible cylindrical-end bundles $\widehat{P}$.  If $\alpha$ is generator for the chain group defining $I^w_*(K_0)$, let
\begin{displaymath}
M_w(\widehat{W}_K)(\alpha), \quad M_w(\widehat{W}'_K)(\alpha)
\end{displaymath}
denote the $L^2$-ASD-moduli spaces with limit $\alpha$ at the cylindrical-ends. (We refer to \cite{D} for details of these constructions.) Without loss of generality we assume these are regular given appropriate perturbations of the equations and suppress the perturbations from our notation. At the chain level, we have
\begin{displaymath}
\delta_K(\alpha) = \# M^{(0)}_w(\widehat{W}_K)(\alpha), \quad 
\delta'_K = \sum_i\# M^{(0)}_w(\widehat{W}'_K)(\alpha_i)\alpha_i
\end{displaymath}
where the super script denotes the 0-dimensional moduli space. Following  arguments in \cite{Fy} (also \cite{D}), it is seen that 
$\delta_K(\partial \alpha)=0$ and $\partial \delta'_K=0$ where $\partial$ is the boundary operator on the Floer chain complex. Therefore these objects are defined at the level of homology.

Since we are in the context of non-trivial $SO(3)$-bundles on $K_0$, the action of the point class via the $\mu$-map is well-defined and determines a degree -4 map on $I^w_*(K_0)$. We denote this action by $u$. (Remark: in comparing our $u$-map with that in \cite{Fy}  the $u$ map there equals $4\mu(\mbox{pt})$.)
Note that $\delta_K$ lies in the dual space of $I^w_*(K_0)$, so the action of $u$ on $\delta_K$ is by pull-back. 

Instanton homology satisfies the non-degenerate pairing
\begin{equation}\label{floer-duality}
\langle ,  \rangle: I^w_*(Y) \times I^w_{k-*}(\overline{Y})\to \mathbb{Q}.
\end{equation}
Roughly speaking, reversing the orientation of $Y$ reverses the direction of the Chern-Simons gradient flow and thus the direction of the differentials in the chain complex for $I^w_*(Y)$. 
The canonical pairing $\langle x_i,x_j\rangle=\delta_{i,j}$, the Kroneckor delta function,  on a set of generators for the chain groups defines the pairing. The pairing gives us the relation
\begin{displaymath}
\langle \delta'_K, \alpha\rangle = \delta_{K^*}(\alpha),
\end{displaymath}
remembering that there is a change of orientation when we switch the role from an ingoing to outgoing boundary component.
The above pairing 
allows us to identify the pull-back action as the same as the $u$-map on $I^w_*(\overline{K}_0) = I^w_*(K_0)^*$. 

The Floer cohomology $I^*_w(Y)$ is defined by considering the dual of the complex
that defines $I^w_*(Y)$. Via the pairing above we have 
\begin{displaymath}
I^*_w(Y) \cong I^w_*(Y)^*\cong I^w_{k-*}(\overline{Y}).
\end{displaymath}
We may simply define the cohomology as the last group on the  right. 

\begin{defn} Let 
\begin{eqnarray*}
{\varphi}^+(K) &=& \mbox{\rm dim span}_{\mathbb{Q}}
\{{\delta}_K, u^2{\delta}_K, u^4{\delta}_K, \dots\}\\
{{\varphi}^-}(K) &=& \mbox{\rm dim span}_{\mathbb{Q}}
\{{\delta}'_K, u^2{\delta}'_K, u^4{\delta}'_K, \dots\}\\
{\varphi}(K) &=& {\varphi}^+(K^*) = {\varphi}^-(K).
\end{eqnarray*}
\end{defn}
We shall work primarily with ${\varphi}^+$ and recover ${\varphi}$ by passing to the mirror knot.

\begin{remark}
Strictly speaking we should orient all our moduli spaces to define the classes $\delta_K$ and $\delta'_K$, otherwise these are only defined to an overall $\pm$-sign. Since all the results in this paper are not dependent on the actual sign we shall proceed without doing so at no loss.
\end{remark}

\section{Proof of Theorem~\ref{thrm-one}}\label{proofthrm-one}

\subsection{Immersed negative concordance}

Let $K$ and $K'$ be two (oriented) knots in $S^3$. Recall that a concordance between $K$ and $K'$ is the existence of an (oriented) embedded annulus $\sigma\cong [0,1]\times S^1$ in $[0,1]\times S^3$ such that $\partial \sigma=\{1\}\times K' -\{0\}\times K$. We weaken the definition of concordance between $K$ and $K'$ by allowing $\sigma$ to be immersed with only negative self-intersection points. Let us write this relation as $K\ge K'$. This relation is not symmetric, however $K\ge K'$ implies ${K'}^*\ge K^*$. If both $K\ge K'$ and $K'\ge K$ then $\sigma$ must have no self-intersection points and $K$ and $K'$ are concordant. If $K'$ is obtained from $K$ by changing a positive crossing in a projection, then $K\ge K'$ with $\sigma$ possessing a single  transverse negative self intersection point.   

\begin{remark} Our convention is that the orientation on a product $[0,1]\times Y$ where $Y$ is an oriented 3-manifold, is $e\wedge\omega$ where $e$ is the standard orientation on $(0,1)$ and $\omega$ is the given orientation on $Y$. Thus if $X$ is an oriented 4-manifold with boundary, a collar of the boundary is orientation preserving  diffeomorphic to $(0,1]\times \partial X$.
\end{remark}

\subsection{Proof of Theorem~\ref{thrm-one}(\ref{thrm-item1})}

Let us now write $\sigma:K\to K'$ to denote the weaker notion of concordance where $\sigma$ is the immersed annulus allowed to have negative transverse self-intersection points, as above. Then $\sigma$ will always induce a map
\begin{equation}\label{induced}
Z_{\sigma *}: I^w_*(K_0)\to I^{w}_*(K'_0)
\end{equation}
by resolving the negative self intersection points in $\sigma$. We recall the procedure of removing the negative self intersection points. Cut out small 4-balls around the self-intersection points to give $\sigma_0$ and replacing each 4-ball by the complement of a 4-ball in $\overline{\mathbb{CP}}^2$ in the following way. 
In $\overline{\mathbb{CP}}^2$ the essential 2-sphere $S$ when pushed off itself to $S'$ intersects $S$ at a single point $p$ with a $-1$-intersection point, so $S\cdot(-S')=+1$. Take the complement of the 4-ball around $p$, it now seen that $\sigma_0$ can be closed off by $S$ and $-S'$ to give the resolved surface $\widetilde{\sigma}$ that algebraically meets the essential sphere in $\overline{\mathbb{CP}}^2$ zero times. 
So $\widetilde{\sigma}$ is an embedded surface in $([0,1]\times S^3)\#n\overline{\mathbb{CP}}^2$, where $n$ is the number of self-intersection points.

Attach a 2-handle to $H$ along $K$ with zero framing where we recall $K$ is the boundary of $\widetilde{\sigma}$ in $\{0\}\times S^3$. Thus we may close off $\widetilde{\sigma}$ in the 2-handle to give a surface $\widehat{\sigma}$. By removing a tubular neighborhood of $\widehat{\sigma}$ in $H$ we obtain a cobordism
\[
Z_{\sigma}: K_0 \to K'_0.
\]
There is a 2-disk in the attached handle that meets the core of the handle transversely once. Via Poincare duality this disk defines a $w_2$ class in $Z_{\sigma}$ that restricts to the unique non-trivial $w_2$ on $K_0$ and $K'_0$. In this way we see there is an induced map as stated in (\ref{induced}).

If we were to close off the $\{1\}\times S^3$ boundary of $H$ then we obtain $W_K\#n\overline{\mathbb{CP}}^2$ where $W_K$ is the defining cobordism for $\delta_K$.  Then 
\[
W_K\#n\overline{\mathbb{CP}}^2 = W_{K'}\circ Z_{\sigma}.
\]
Let $\overline{\delta}_K$ denote the zero-dimensional Donaldson invariant for $W_K\#n\overline{\mathbb{CP}}^2$. Then the blow-up relations for Donaldson invariants tell us that $\overline{\delta}_K=\delta_K$ since the $w_2$-class is zero on the $\overline{\mathbb{CP}}^2$ factors (see for example \cite{FM}).
Then it follows that
\[
Z_{\sigma}^*(\delta_{K'}) = \delta_K.
\]
Since $Z_{\sigma}^*$ commutes with the $u$-map we see that ${\varphi}^+(K')\leq {\varphi}^+(K)$. Passing to the mirror knots we find $K^*\leq K'^*$ and ${\varphi}(K^*)\leq {\varphi}(K'^*)$.  
If $K$ is concordant to $K'$ then both $K\ge K'$ and $K'\ge K$ so ${\varphi}(K)={\varphi}(K')$.

\subsection{Some results on the instanton homology for $S^1\times\Sigma$}

In order to prove the 4-ball genus we need to recall some results of Munoz \cite{M} that will be the basis of the genus bound. Let $\Sigma$ be a surface of genus $g$. Let $w=w_2$ be chosen to be such that $w$ is the Poincare dual to the 2-dimensional homology class represented by 
$\{\mbox{pt}\}\times\Sigma\subset S^1\times\Sigma$. In addition to the $u$-map acting on ${I}^w_*(S^1\times\Sigma)$ we have actions defined by $\mu(\Sigma)$ and $\mu(\xi_i)$ where $\xi_i$, $i=1,\dots, 2g$ is a symplectic basis $\xi_i\cdot\xi_{g+j}=\delta_{i,j}$ for $H_1(\Sigma)$. (We  ignore the action by $\mu(S^1)$ because it is the zero map.) 
Denote $\mu(\Sigma)$ by $\sigma$ and $\mu(\xi_i)$ by $\gamma_i$. (In \cite{M} $\sigma$ is denoted by $\alpha$ and the $u$-map by $\beta$ modulo normalizations.)

Let $Z=D^2\times\Sigma$. Since $Z:\emptyset\to D^2\times\Sigma$ this determines an element $\psi$ in $I^w_*(D^2\times\Sigma)$ by considering the zero-dimensional Donaldson invariant over $\widehat{Z}$. If we were to instead let $Z:D^2\times\Sigma\to\emptyset$ we have an element $\psi^*:{I}^w_*(S^1\times\Sigma)\to\mathbb{Q}$. (Strictly $\psi$ is dual to an element $\overline{\psi}^*:{I}^w_*(\overline{S^1\times\Sigma})\to\mathbb{Q}$.) Let us now pass to the $\mathbb{Z}/4$-graded groups $\widetilde{I}^w_*(S^1\times\Sigma)$ by taking the quotient by the degree 4 involution represented by the Poincare dual of $\{\mbox{pt}\}\times\Sigma$, as a class in $H^1(S^1\times\Sigma, \mathbb{Z}/2)$.

Let $\mathbb{A}=\mathbb{A}[\sigma, u, \gamma_i]$ be the graded commutative ring over $\mathbb{Q}$ generated by $\sigma$, $u$ and $\gamma_i$ where $\sigma$ is degree 2, $u$ degree 4 and $\gamma_i$ degree 3. We can then consider $\widetilde{I}^w_*(S^1\times\Sigma)$ as a module over the ring $\mathbb{A}$.
\cite[p.520]{M} states that the elements
$\widetilde{\psi},\sigma\widetilde{\psi}, u\widetilde{\psi}, \gamma_i\widetilde{\psi}, i=1,\dots, 2g$
are linearly independent and generate $\widetilde{I}^*_w(S^1\times\Sigma)$ as a module over $\mathbb{A}$. Therefore if we set $\widetilde{\psi}$ to be the identity element, this will make $\widetilde{I}^*_w(S^1\times\Sigma)$ into a ring generated by $\sigma$, $u$ and the $\gamma_i$.

Let $\gamma$ be the element  $\sum_i\gamma_i\gamma_{g+i}$ of degree 6 in $\mathbb{A}$. Then $\mathbb{Q}[\sigma, u,\gamma]$ is a commutative subring of $\mathbb{A}$. The submodule $\mathbb{Q}[\sigma, u,\gamma]\widetilde{\psi}$ of $\widetilde{I}^w_*(S^1\times\Sigma)$ corresponds to the subspace that is invariant under diffeomorphisms  induced by the mapping class group of $\Sigma$. \cite{M} determines the ideal of relations $J_g$ for this submodule i.e.
\begin{displaymath}
\mathbb{Q}[\sigma, u,\gamma]\widetilde{\psi} \cong\mathbb{Q}[\sigma, u,\gamma]/J_g.
\end{displaymath}
Without going into the complete analysis in, for our purposes we only need the following relation \cite[proof of Prop.~20]{M}:

\begin{prop}\label{prop-one}
There is a polynomial $p(\sigma, u,\gamma)$ such that
\[
(\Pi_{i=1}^g (u +(-1)^i4) - \gamma\cdot p(\sigma,u,\gamma))\widetilde{\psi} = 0
\]
in $\mathbb{Q}[\sigma, u,\gamma]\widetilde{\psi}\subset \widetilde{I}^w_*(S^1\times\Sigma)$. Therefore in 
$I^w_*(S^1\times\Sigma)$ the following relation holds:
\[
(u^2-4)^{\lceil g/2\rceil}\psi -\gamma q(\sigma, u, \gamma)\psi =  0
\]
for some $q(\sigma, u, \gamma)$.
\end{prop}
We remark that for $g=1, 2$ the polynomial $p(\sigma, u,\gamma)=0$. We also note that \cite{M} works over $\mathbb{C}$ but the results mentioned  above hold over $\mathbb{Q}$. We have different constants from \cite{M} in the proposition because of different normalizations for the $u$-map.

It will be convenient to work dually, and consider $\psi^*$ as above. The non-degenerate pairing (\ref{floer-duality})
implies that the 2nd part of Proposition~\ref{prop-one} also holds for $\psi^*$.

\subsection{Proof of Theorem~\ref{thrm-one}(\ref{thrm-item2})}

We need a small remark regarding of the action of the $\mu$-maps and cobordisms and 1-dimensional classes. A 1-dimensional homology class $\xi\in H_1(Y,\mathbb{Z})$ acts with degree $-3$ on $I^w_*(Y)$ through the action of the $\mu$-map.
Suppose we have a cobordism $(W, \bar{w}): (Y, w)\to(Y',w')$ by which we mean a cobordism with 
$\bar{w}|_{Y}=w$ and $\bar{w}|_{Y'}=w'$. Assume $\xi$ is a 1-dimensional homology class in $Y$ and $\xi'$ a 1-dimensional homology class in $Y'$ such that $\xi\sim \xi'$ in $W$. Then 
$W_*$ is compatible with the action of $\xi, \xi'$ in the following sense
\begin{displaymath}
W_*(\xi \alpha) = \xi' W_*(\alpha).
\end{displaymath}

Let $W_K:K_0\to\emptyset$ be the defining cobordism for $\delta_K$. Assume $K$ bounds an immersed surface $\Sigma_0$ of genus g in the 4-ball with $\Sigma$ having transverse self-intersections points that are negative. As before, we may remove the negative self-intersection points by connect sums with $\overline{\mathbb{CP}}^2$. Since this makes no essential difference to the argument, without loss we can assume $\Sigma_0$ is embedded. $\Sigma_0$ is closed off inside $W_K$ to give a closed surface $\Sigma$ of genus g. (It may be easier to see this by considering instead $W'_{K^*}$ which is a 4-ball with a 0-framed 2-handle attached along $K^*$. $\Sigma$ corresponds to closing off $\Sigma_0$ with the core of the 2-handle. As mentioned previously $W'_{K^*}$ is the same as $W_K$ but with ingoing and outgoing boundaries switched.) We can identify a tubular neighborhood of $\Sigma$ in $W_K$ with $H_g:\partial(D^2\times\Sigma)\to\emptyset$ and we have a factorization
\begin{displaymath}
W_K = H_g\circ Z, \quad Z: K_0\to S^1\times\Sigma.
\end{displaymath}
Consequently,
$Z^*\psi^* = \delta_K$.
From Prop.~\ref{prop-one} 
\begin{displaymath}
0=Z^*((u^2-4)^{\lceil g/2\rceil}\psi^* -\gamma q(\sigma, u, \gamma)\psi^* ) = (u^2-4)^{\lceil g/2\rceil}\delta_K,
\end{displaymath}
where we use the fact that $\gamma$ being constructed from the 1-dimensional classes $\xi_i$ in 
$\{\mbox{pt}\}\times\Sigma$ must be homologically trivial in $Z$. This can be see by from the long exact sequence for the pair $(W_K, W_K - \Sigma)$ and excision. 
The condition $(u^2-4)^{\lceil g/2\rceil}\delta_K=0$ implies 
\begin{displaymath}
{\varphi}^+(K) =
\mbox{dim span}_{\mathbb{Q}}
\{\delta_K, u^2\delta_K, u^4\delta_K, \dots\} \leq \lceil g/2\rceil.
\end{displaymath}
By passing to the mirror $K^*$,  the genus bound proof is complete.

\subsection{Proof of Theorem~\ref{thrm-one}(\ref{thrm-item3})} 

Given a homologically trivial knot $K$ in $Y$ there is an exact sequence due to Floer:
\begin{displaymath}
\longrightarrow
{I}_i^{w}(Y)
\stackrel{a}{\longrightarrow}
{I}_i^w(K_{-1})
\stackrel{b}{\longrightarrow}
{I}_i^{w+w'}(K_0)
\stackrel{c}{\longrightarrow}
{I}_{i-1}^w(Y)
\longrightarrow
\end{displaymath}
where the maps are induced by the various surgery cobordisms. These are obtained by adding a 2-handle to the product cobordism  at each stage \cite{BD}. The maps $a$, $b$ are degree zero and $c$ is degree $-1$ by convention, i.e. choosing the grading on $I^{w+w'}_*(K_0)$ to be compatible. $w'$ is the class unique class introduced by zero-surgery, we can geometrically realize it as a small circle linking $K$. We remark that this additional $w'$ is necessary in the application of Floer's exact sequence, even as $w$ is non-trivial in $K_0$.

Let $K_+$ and $K_{-}$ be two knots that only differ in a crossing change. We can think of this crossing change as being effected by $-1$-surgery along a small circle $C$ that bounds a disk $D$. This disk meets $K_+$ near the crossing change geometrically in exactly two points and algebraically zero times (after an orientation of $K_+$ is given). Then we can think of $K_{-}$ as being the knot $K_+$ in $C_{-1}$. Note that $C$ is homologically trivial in the complement of $K_+$: the disk $D$ meets $K_+$ in two points and after removing small neighborhoods of the intersections points we can tube the disc around one half of $K_+$ to get a genus one surface $T_0$ that bounds $C$ (here we have used the property that the algebraic intersection of $D$ with $K_+$ is zero). In 
$(K_{+0};C_0)$ this genus one surface $T_0$ is closed off to give an embedded surface $T$.

Applying Floer's exact triangle to the knot $C$ we obtain
\begin{displaymath}
\longrightarrow
{I}_i^{w}(K_{+0})
\stackrel{a}{\longrightarrow}
{I}_i^w(K_{+0};C_{-1})
\stackrel{b}{\longrightarrow}
{I}_i^{w+w'}(K_{+0}; C_0)
\stackrel{c}{\longrightarrow}
{I}_{i-1}^w(K_{+0})
\longrightarrow
\end{displaymath}
where $w'$ is the class that evaluates non-trivially over $T$. Observe that $(K_{+0};C_{-1}) \cong K_{-0}$. 

\begin{lem} 
The crossing change cobordism induced by surgery along $C$ above satisfies $a^*\delta_{K_{-}} = \delta_{K_{+}}$.
\end{lem}

Assume the lemma. It immediately follows that ${\varphi}^+(K_+)\leq {\varphi}^+(K_{-})$.
As  mentioned above $(K_{0+},C_0)$ contains an embedded torus $T$ for which $w'$ evaluates non-trivially and odd over $T$. This forces
$(u^2-4)=0$ in ${I}^{w+w'}(K_{+0};C_0)$. For a proof we refer to \cite[Lemma~4]{Fy}. 
Suppose there is a polynomial $p_n$ of degree $n$ such that $p_n(u^2)\delta_{K_{+}}=0$ in $I^w_*(K_{+0})$. Then $a^* (p_n(u^2)\delta_{K_{-}})=p_n(u^2)a^*\delta_{K_{-}}=0$ by the above lemma. Since we are working over $\mathbb{Q}$ the dual sequence is also exact and we see that $p_n(u^2)\delta_{K_{-}}=b^*(z)$ for some $z$. 
Then $(u^2-4)p_n(u)\delta_{K_-}=0$  and therefore ${\varphi}^+(K_{-})\leq {\varphi}^+(K_{+})+1$. Passing to the mirror $K^*$ the crossing change bound proof follows.

\proofoflem 
Let $W_{K_+}$ be the defining cobordism for $\delta_{K_+}$ and $W^{\circ}_{K_+}:K_{+0}\to S^3$ be obtained from
$W_{K_+}$ by removing an open 4-ball, as in (\ref{cob-one}). Recall $W^{\circ}_{K_+}$ is 
obtained by attaching a two handle to $K_{+0}\times[0,1]$ along $K'\times\{1\}$ where $K'$ is the knot described preceding (\ref{cob-one}). 
If we follow this by a 2-handle attached along $C=C\times\{1\}$ with a $-1$ framing we will get 
$\widetilde{W}^{\circ}:K_{+0}\to S^3$ with 
$\widetilde{W}^{\circ}\cong W^{\circ}_{K_+}\#\overline{\mathbb{CP}}^2$. 
(The effect of $W^{\circ}_{K_+}$ is to erase $K_{+}$ itself i.e. $\infty$-surgery, so $C$ will be the unknot in $S^3$.)  Fill in the $S^3$ with a 4-ball for $\widetilde{W}^{\circ}$ to give  
$\widetilde{W}\cong W_{K_{+}}\#\overline{\mathbb{CP}}^2$. Let $\delta_{\widetilde{W}}:I^w_*(K_{+0})\to\mathbb{Z}$ be the 0-dimensional Donaldson invariant for $\widetilde{W}$ where $w_2$ is non-trivial on $K_{+0}$ but trivial on $\overline{\mathbb{CP}}^2$; then $\delta_{K_{+}}=\delta_{\widetilde{W}}$.

On the other hand $\widetilde{W}$ can also be factorized as $\widetilde{W}=W_{K_{-}}\circ H$ where $W_{K_{-}}$ is the defining cobordism for $\delta_{K_{-}}$ and $H$ the cobordism that induces $a$ in the exact sequence. This  is obtained by switching the order of attaching the 2-handles. It now follows that 
$a^*\delta_{K_{-}}=\delta_{\widetilde{W}}=\delta_{K_{+}}$ and this proves the lemma.

\section{Proof of Theorem~\ref{thrm-one}(\ref{thrm-item4})}

\subsection{A cobordism map}
We refer to section~\ref{subsection-connect} for the discussion of the connect sum complex for instanton homology, which we shall assume. Let $(Y,w)$ and $(Y',w')$ be \textit{non-trivial admissible} pairs i.e. where $w, w'$ evaluate non-trivially on some closed oriented surface in $Y, Y'$ respectively. Let $(C,\partial)$, $(C', \partial')$ denote the Floer complex for $I^w_*(Y)$, $I^{w'}_*(Y')$ respectively. Then the Floer complex for $I^{w+w'}_*(Y\sqcup Y')$ is given by 
\[
\widetilde{C}=(C\otimes C')_*\oplus (C\otimes C')[3]_*,\quad  \widetilde{\partial}=\partial \otimes I +\epsilon I\otimes \partial' + \Phi
\]
where
\[
\Phi: (C\otimes C')_*\to (C\otimes C')[3]_*, \quad \Phi= u\otimes I - I\otimes u.
\]
This gives 
\[
I^{w+w'}_*(Y\sqcup Y') = \text{ker}\widehat{\Phi}\oplus \text{coker}\widehat{\Phi}
\]
where
\[
\widehat{\Phi}: I^w(Y)\otimes I^{w'}(Y') \to I^w(Y)\otimes I^{w'}(Y')[3], \quad\widehat{\Phi}=u\otimes I-I\otimes u.
\]
Assume that $(Y'',w'')$ is also admissible and we have a connected cobordism
\begin{equation}\label{W-cob-a}
W:Y\sqcup Y'\to Y''.
\end{equation}
We allow the possibility that $Y''$ is empty.  We assume we have a $w_2$ class $\overline{w}$ that restricts to classes $w$, $w'$, $w''$ on $Y$, $Y'$, $Y''$ respectively.   
Following the discussion in Section~\ref{subsection-connect} we have a map
\[
W_*: I^{w+w'}_*(Y\sqcup Y')\to I^{w''}_*(Y'').
\]
If $Y'' =\emptyset$ interpret $I^{w''}_*(Y'')$ as $\mathbb{Q}$. In the case where $W$ is the standard cobordism $X: Y\sqcup Y'\to Y\# Y'$, this map is an isomorphism by Fukaya's result \cite{Fu} (see Theorem~\ref{connectsum} in section~\ref{subsection-connect} below).

The action of the $u$-map extends to $I^{w+w'}_*(Y\sqcup Y')$ when we define
\[
u(a\otimes b, c\otimes d) =(ua\otimes b, uc\otimes d)
\]
for
$a\otimes b \in C(Y)\otimes C(Y')$ and $c\otimes d \in C(Y)\otimes C(Y'))[3]$, as this action commutes with the boundary operator of the complex. 
Observe  that $\widehat{\Phi}$ gives the additional symmetry relations
\[
(ua\otimes b, uc\otimes d) = (a\otimes ub, uc\otimes d) = (a\otimes ub, c\otimes ud)
\]
\textit{provided} $(a\otimes b, c\otimes d)$ is in $I^{w+w'}_*(Y\sqcup Y')$. (Remark: in general we do not have $(ua\otimes b,0)= (a\otimes ub,0)$. This  holds if and only if $\widehat{\Phi}(a\otimes b)=0$ i.e. $a\otimes b$ is a cycle in the complex.)

\begin{lem} \label{W-lemma}
$W_*$ commutes with the $u$-map. We may allow $Y''$ to be empty, in which case we interpret $I^{w''}_*(Y'')$ as $\mathbb{Q}$.
\end{lem}

\proof
This is a standard argument involving looking at the ends of a non-compact 1-dimensional moduli space, as done in \cite{Fy} together with some results on $\mu(\widehat{\gamma})$ in \cite{D}. That $W_*$ defines a chain map is in \cite{Sca}. To show $W_*$ commutes with the $u$-map we follow similar arguments and we only outline the key points.  
At the chain level we have $W_*= W^{\gamma}_*+ W^0_*$ with $W^{\gamma}_*: C(Y)\otimes C(Y')\to C(Y'')$ and $W^0_*: C(Y)\oplus C(Y')[3]\to C(Y'')$. These are defined by the following rules
\begin{eqnarray*}
&&\langle W^{\gamma}_*(a\otimes b ),c\rangle = \#M^{(3)}(\widehat{W})(a, b;c)\cap\mu(\widehat{\gamma})\\
&&\langle W^0_*(a\otimes b), c\rangle = \#M^{(0)}(\widehat{W})(a, b;c).
\end{eqnarray*}
In the case where $Y''=\emptyset$ we can, by removing a 4-ball replace $Y''$ with $S^3$ and let $c$ be the trivial connection $\theta$; then $C(Y'')$ will be replaced by $\mathbb{Q}\langle\theta\rangle$.
Recall that $\gamma$ is an embedded curve that connects $Y$ with $Y'$ in $W$; $\widehat{\gamma}$ is $\gamma$ (after possibly a reparameterized) extended over the cylindrical-end manifold $\widehat{W}$.
To demonstrate that $W^0_*$ commutes with the $u$-action (this is essentially proven in \cite{Fy}), let $r:\mathbb{R}\to \widehat{W}$ such that for $s<-1$, $r(s)=(s, \mbox{pt})$  on the cylindrical-end $(-\infty, 0]\times Y$, and for $s>1$, $r(s)=(s, \mbox{pt})$ on the cylindrical-end $[0,\infty)\times Y''$. Consider the 1-dimensional moduli space
\[
\bigcup_s \{s\}\times M^{(4)}_{\overline{w}}(\widehat{W})(a, b;c)\cap \mu(r(s))\subset 
\mathbb{R}\times M^{(4)}_{\overline{w}}(\widehat{W})(a, b;c).
\]
Non-compactness occurs when instantons bubble-off at $t=-\infty$ at the ingoing cylindrical-ends $(-\infty,0]\times \partial_{-} W$ of $\widehat{W}$ and also bubbling off at the $t=+\infty$ at  the outgoing cylindrical-end $[0,+\infty)\times\partial_+W$ of $\widehat{W}$. The ends are orientation preserving diffeomorphic to
\begin{alignat*}{2}
&\mathbb{R}\times M^{(4)}_{w}(\mathbb{R}\times Y)(a,\tau)\cap\mu(\text{pt})\times
M^{(0)}_{\overline{w}}(\widehat{W})(\tau, b;c),&\quad (s, t)\to (-\infty, -\infty)\\
& M^{(0)}_{\overline{w}}(\widehat{W})(a, b;\tau)\times\mathbb{R}\times M^{(4)}_{w}(\mathbb{R}\times Y)(\tau,c)\cap\mu(\text{pt}), &\quad (s,t)\to (+\infty,+\infty).
\end{alignat*}
This gives the relation $\langle W^0_*(u a\otimes b),c\rangle -\langle uW^0_*(a\otimes b),c\rangle=0$. 
To show $W^{\gamma}_*$ commutes with the $u$-map consider the 1-dimensional moduli space
\[
\bigcup_s \{s\}\times M^{(7)}_{\overline{w}}(\widehat{W})(a, b;c)\cap\mu(r(s))\cap\mu(\widehat{\gamma})\subset
\mathbb{R}\times M^{(7)}_{\overline{w}}(\widehat{W})(a, b;c).
\]
Let $V$ be $Y$ or $Y'$  and 
\[
\widehat{M}^{(3)}_w(\mathbb{R}\times V)(a, b) = M^{(4)}_w(\mathbb{R}\times V)(a, b)/\mathbb{R}
\]
be the reduced moduli space, i.e.  the action by $\mathbb{R}$ is by translation. Let $\overline{\gamma}(t)$be the path $(t, \text{pt})$ in $\mathbb{R}\times V$. 
One key fact we need from \cite{D} is
\begin{equation}\label{mu-hol}
\#\widehat{M}^{(3)}_w(\mathbb{R}\times V)(a, b) \cap\mu(\overline{\gamma}) = \# M^{(4)}_w(\mathbb{R}\times V)(a, b)\cap\mu(\text{pt}).
\end{equation}
Then the ends of the moduli space are orientation preserving diffeomorphic to
\begin{alignat*}{2}
&\mathbb{R}\times M^{(7)}_w(\mathbb{R}\times Y)(a;\rho)\cap\mu(\text{pt})\cap\mu(\widehat{\gamma}_{+})\times M^{(0)}_{\overline{w}}(\widehat{W})(\rho, b;c), &(s,t)\to(-\infty,-\infty)\\
&\mathbb{R}\times M^{(4)}_w(\mathbb{R}\times Y)(a;\rho)\cap\mu(\text{pt}) \times \widehat{M}^{(3)}_{w'}(\mathbb{R}\times Y')(b,\rho')\cap\mu(\widehat{\gamma}_{-})\times M^{(0)}_{\overline{w}}(\widehat{W})(\rho, \rho';c), &(s,t)\to(-\infty,-\infty)\\
&\mathbb{R}\times M^{(4)}_w(\mathbb{R}\times Y)(a;\rho)\cap\mu(\text{pt}) \times M^{(3)}_{\overline{w}}(\widehat{W})(\rho, b;c)\cap \mu(\widehat{\gamma}), &\quad (s,t)\to(-\infty,-\infty)\\
&\mathbb{R}\times\widehat{M}^{(3)}_w(\mathbb{R}\times Y)(a,\rho)\cap\mu(\widehat{\gamma}_{+})\times M^{(0)}_{\overline{w}}(\widehat{W})(\rho, b;\rho')\times M^{(4)}_{w''}(\mathbb{R}\times Y'')(\rho',c)\cap\mu(\text{pt}) , &(s,t)\to(+\infty,\pm\infty)\\
&\mathbb{R}\times\widehat{M}^{(3)}_{w'}(\mathbb{R}\times Y')(b,\rho)\cap\mu(\widehat{\gamma}_{-})\times M^{(0)}_{\overline{w}}(\widehat{W})(a,\rho;\rho')\times M^{(4)}_{w''}(\mathbb{R}\times Y'')(\rho',c)\cap\mu(\text{pt}) , &(s,t)\to(+\infty,\pm\infty)\\
&\mathbb{R}\times M^{(3)}_{\overline{w}}(\widehat{W})(a, b;\rho)\cap \mu(\widehat{\gamma})\times M^{(4)}_{w''}(\mathbb{R}\times Y'')(\rho',c)\cap\mu(\text{pt}) , &(s,t)\to(+\infty,\pm\infty).
\end{alignat*}
Here $\widehat{\gamma}_{+}$ is a path $(t, \text{pt})$ in $\mathbb{R}\times Y$ and 
$\widehat{\gamma}_{-}$ is a path $(-t, \text{pt})$ in $\mathbb{R}\times Y'$.  
The ends give the relation 
\begin{multline}\label{ends-hol-u}
\qquad 0=\langle W_*(u^2a\otimes b),c\rangle -\langle W_*(ua\otimes ub),c\rangle 
+\langle W_*^{\gamma}(u a\otimes b), c\rangle\\
 +\langle uW_*(ua\otimes b), c\rangle
-\langle uW_*(a\otimes ub), c\rangle -\langle uW^{\gamma}_*(a\otimes b),c\rangle.
\end{multline}
If $\widehat{\Phi}(a\otimes b)=0$ then  $ua\otimes b-a\otimes ub=0$ and we see
the 4th and 5th terms on the right-side sum to zero. Likewise since the action by $u$ is a chain map, $ua\otimes b$ is also a cycle, so the first and 2nd terms sum to zero as well. This leaves 
\[
W_*^{\gamma}(u a\otimes b) = uW^{\gamma}_*(a\otimes b)
\]
and the lemma follows.

\subsection{Completion of proof}

\begin{lem}\label{main-lem} 
Let $W_K$, $W_{K'}$ be the defining cobordisms for $\delta_K$, $\delta_{K'}$ respectively. Denote by $W_{K, K'}$ the boundary connected sum of $W_K$ and $W_{K'}$, and $\delta_{K,K'}$ the corresponding zero dimensional Donaldson invariant. Let $(u^2-4)^n=0$ in both $I^w_*(K_0)$ and $I^{w'}_*(K'_0)$ for some $n\ge 1$. Assume $ {\varphi}(K)+{\varphi}(K')\ge n+1$. Then $(u^2-4)^l\delta_{K,K'}\neq 0$ where $l={\varphi}(K)+{\varphi}(K')-n-1$.
\end{lem}

\proof Rather than work directly with $W_{K,K'}$ we compose with the cobordism $X: K_0\sqcup K'_0\to K_0\# K'_0$. Then it is seen that
\begin{equation}\label{W-connectsum}
W_{K,K'}\circ X \cong W_K \# W_{K'}.
\end{equation}
$X$ determines an isomorphism $X_*:I^{w+w'}_*(K_0\sqcup K'_0)$ with $I^{w+w'}_*(K_0\#K'_0)$ via Fukaya's Theorem. We have 
\[
I^{w+w'}_*(K_0\sqcup K'_0)=\text{ker}\widehat{\Phi}\oplus \text{coker}\widehat{\Phi}
\]
where
\[
\widehat{\Phi}: I^w(K_0)\otimes I^{w'}(K'_0) \to I^w(K_0)\otimes I^{w'}(K'_0)[3], \quad\widehat{\Phi}=u\otimes I-I\otimes u.
\]
Having identified $I^{w+w'}_*(K_0\sqcup K'_0)$ with $I^{w+w'}_*(K_0\#K'_0)$ via $X_*$ we can now regard $\delta_{K,K'}$ as a map $I^{w+w'}_*(K_0\sqcup K'_0)\to\mathbb{Q}$.

\begin{claim} Let $\alpha=\sum_i a_i\otimes b_i$ be a class in $I^w_*(K_0)\otimes I^{w'}_*(K'_0)$ such that $\widehat{\Phi}(\alpha)=0$. Then 
\[
\delta_{K,K'}(\alpha) = \frac{1}{2}\sum_i\delta_K(a_i)\delta_{K'}(b_i).
\]
\end{claim}
To see the claim: at the chain level 
\[
\delta_{K,K'}(a\otimes b)=\# M^{(3)}_w(\widehat{W}_K\#\widehat{W}_{K'})(a, b)\cap\mu(\widehat{\gamma})
\]
where ${\gamma}$ is an embedded path that connects $K_0$ with $K'_0$ inside $W_K \# W_{K'}$. 
In the connected sum region we may consider a long cylindrical tube of length $T>0$ modeled on $[-T,T]\times S^3$ with a standard metric on $S^3$. Let $T\gg 1$. By standard gluing arguments for such $T$ the moduli space is orientation preserving diffeomorphic to 
\[
M^{(0)}_w(\widehat{W}_K)(a)\times M^{(0)}_{w'}(\widehat{W}_{K'})(b)\times SO(3)
\]
where $SO(3)$ represents the gluing parameter. Now we use the crucial fact that $\mu(\widehat{\gamma})$ evaluates to be $1/2$ on the gluing parameter (see \cite{D}). The claim now follows.

Our assumption is the existence of $n\ge 1$ such that $(u^2-4)^n=0$ for $I^w_*(K_0)$ and $I^{w'}_*(K'_0)$. 
Thus we have the following identity. For convenience let $u_1=u\otimes I$ and $u_2=I\otimes u$.
\[
0=(u_1^2-4)^n-(u_2^2-4)^n = (u_1 - u_2)(u_1+u_2)\sum_{i=0}^{n-1} (u_1^2-4)^i(u_2^2-4)^{n-i}.
\]
The condition that ${\varphi}(K)=k>0$ can be equivalently stated as $(u^2-4)^i\delta_K\neq 0$ for $0\le i\le k-1$ and $(u^2-4)^i\delta_K=0$ for $i\ge k$. Since the vectors $\delta_K,(u^2-4)\delta_K,\dots,(u^2-4)^{k-1}\delta_K$ are independent, there must exist an $a$ such that $\delta_K((u^2-4)^i a)=0$ for $0\le i\le k-2$ and $\delta_K((u^2-4)^{k-1} a)\neq 0$.
An identical statement holds for ${\varphi}(K')=k'>0$ and element $b$. 

The $u$ map is an isomorphism and therefore we can write $a=ua'$ for some $a'$. It now follows from the above identity that 
\[
\alpha=\sum_{i=0}^{n-1} (u^2-4)^i ua'\otimes (u^2-4)^{n-1-i}b+(u^2-4)^i a'\otimes (u^2-4)^{n-1-i}ub
\]
satisfies $\widehat{\Phi}(\alpha)=0$ i.e. $\alpha$ represents a class in $I^{w+w'}_*(K_0\sqcup K'_0)$. 
Then by the claim we have
\[
(u^2-4)^l\delta_{K,K'}(\alpha) = \delta_K((u^2-4)^{k-1} a)\delta_{K'}((u^2-4)^{n-1-(k-1)+l}b)\neq 0
\]
when $n-1-(k-1)+l=k'-1$. This proves  the lemma. \qed

\medskip

We may iterate the construction above to three knots. 
Let
\[
W_{K,K',K''} = W_{K}\#_{\partial}W_{K'}\#_{\partial}W_{K''}
\]
and denote the $w_2$ classes by $w$, $w'$ and $w''$. Let $(u^2-4)^n=0$ in $I^w_*(K_0)$, $I^{w'}(K'_0)$ and $I^{w''}_*(K''_0)$ for some $n\ge 1$. We can extend the proof of  Lemma~\ref {main-lem} to this context:  if $k+k'+k''\ge 2n+1$ then 
\[
(u^2-4)^l\delta_{K, K',K''}\neq 0, \quad l=k+k'+k''-2n-1.
\]
We outline the key points. Let $X':(K_0\# K'_0)\sqcup K''_0 \to K_0\#K'_0\# K''_0$ be the standard cobordism just as in the proof of the lemma. Now compose with $X: K_0\sqcup K'_0\to K_0\# K'_0$ to give the cobordism
\[
X'\circ X: K_0\sqcup K'_0\sqcup K''_0\to K_0\#K'_0\# K''_0.
\]
The relevant chain complex is the direct sum of the following four complexes
\begin{align*}
&I^w(K_0)\otimes I^{w'}(K'_0)\otimes I^{w''}(K''_0)\\
&I^w(K_0)\otimes I^{w'}(K'_0)\otimes I^{w''}(K''_0)[3]\\
&I^w(K_0)\otimes I^{w'}(K'_0)\otimes I^{w''}(K''_0)[3]\\
&I^w(K_0)\otimes I^{w'}(K'_0)\otimes I^{w''}(K''_0)[6].
\end{align*}
An element $a\otimes b\otimes c$ is a cycle if it satisfied $(u_1-u_2)(u_1-u_3)(a\otimes b\otimes c) =0$.
By employing the identity
\[
0=[(u^2_1-4)^n - (u^2_2-4)^n][(u^2_1-4)^n-(u^2_3-4)^n]
\]
we find that the element
\[
\alpha'=\sum_{i,j}^{n-1} (u^2_1-4)^{i+j}\otimes (u^2_2-4)^{n-1-i}\otimes (u^2_3-4)^{n-1-j}(u_1+u_2)(u_1+u_3)(a\otimes b\otimes c)
\]
is a cycle. Then 
\[
(u^2-4)^l\delta_{K,K',K''}(\alpha) = \frac{1}{4}\delta_K((u^2-4)^{k-1} a)\delta_{K'}((u^2-4)^{k'-1}b)\delta_{K''}((u^2_3-4)^{k''-1}\neq 0
\]
which happens for $i+j=k-1$, $n-1-i=k'-1$ and $n-1-j+l=k''-1$. This gives $l=k+k'+k''-2n-1$.

Let $K\cup K'\cup K''$ be the link in $S^3$ where each of the $K$, $K'$, $K''$ are contained in disjoint 3-balls.
Then there is an embedded surface $\sigma$ in $[0,1]\times S^3$ such that $\partial\sigma$ is $K\cup K'\cup K''\cup K\# K'\#K''$ where $K \cup K'\cup K''$ is in the ingoing boundary $\{0\}\times S^3$ and $K\# K'\#K''$ is in the outgoing $\{1\}\times S^3$ boundary. Then following the same construction as in the proof of Theorem~\ref{thrm-one}(\ref{thrm-item1}), $\sigma$ determines a factorization
\[
W_{K,K',K''} = W_{K\# K'\#K''}\circ Z_{\sigma}, \quad Z_{\sigma}: K_0\# K'_0\# K''_0\to (K\# K'\# K'')_0.
\]
The $w_2$ class for $(K\# K'\# K'')_0$ is represented by a meridian $m^{\#}$ of $K\# K'\# K''$ and likewise the meridians $m$, $m'$, $m''$ for $K$, $K'$,  $K''$ respectively. In $Z_{\sigma}$,
$m^{\#}$ and $m$ cobounds an annulus, and same goes for $m^{\#}$ and $m'$, and $m^{\#}$ and $m''$. This shows that there is a $w_2$ class on $Z_{\sigma}$ that restricts to $w^{\#}$, $w$, $w'$, $w''$ on the boundaries. (Note that it is at this point that the necessity for the connected sum of three knots comes into play as no $w_2$ with the desired properties will exist if we only consider only two knots.)
Therefore
\[
(u^2-4)^l\delta_{K,K',K''}=Z^*_{\sigma}(u^2-4)^l\delta_{K\#K'\#K''} \neq 0
\]
and Theorem~\ref{thrm-one}(\ref{thrm-item4}) follows.

\section{Results on connected sums}\label{conn}

\subsection{The Poincare homology sphere}

For an integral homology 3-sphere $Y$ denote by $M(\mathbb{R}\times Y)(\alpha,\beta)$ the moduli space of (perturbed) ASD connections that are asymptotic to $\alpha$ at $-\infty$ and $\beta$ at $+\infty$.  
Following \cite{Fy} we have two maps of special importance:
\begin{eqnarray*}
&\delta: C_1(Y)\to\mathbb{Q}, \quad &\delta(\alpha) =\# (M(\mathbb{R}\times Y)(\alpha, \theta))/\mathbb{R}\\
&\delta':\mathbb{Q}\to C_4(Y),\quad &\delta'(1) =\sum_i\# (M(\mathbb{R}\times Y)(\theta,\alpha_i))/\mathbb{R}
\end{eqnarray*}
where the action of $\mathbb{R}$ is by translation, the $\alpha_i$ are representatives that run over a basis for $I_4(Y)$ and $\theta$ is the equivalence class of the trivial connection. These descend to maps
\[
\delta_0: I_1(Y)\to\mathbb{Q}, \quad \delta'_0:\mathbb{Q}\to  I_4(Y).
\]
In the context of $\mathbb{Z}$-homology spheres the $\delta$ and $\delta'$ maps are correction terms to the commutator of the $u$-map and boundary operator at the chain level (\cite{Fy}, \cite{D}) 
\begin{equation}\label{u-chain}
\partial u - u\partial +
\frac{1}{2}\delta'\circ \delta =0.
\end{equation}
This shows the $u$-map does not generally descend to a map at the level of homology. (Note: when comparing formulas our $u$-map is $1/4$ of the $u$-map in \cite{Fy}, 2 times the $U$-map in \cite{D} and $1/2$ of the $v$-map in \cite{Sca}).

Let us now collect from facts about the Poincare sphere from \cite{Fy}. Let $P^{-}$ be the Brieskorn homology sphere $B(2,3,5)$ oriented as the boundary of a 4-manifold with intersection form $-E_8$ i.e. we may write $-E_8:\emptyset\to P^-$. Reversing the ingoing and outgoing component we get $-E^*_8: P^+\to\emptyset$ where $P^+$ now has the reversed orientation. We shall follow this orientation notation convention throughout the paper. 

With our orientation conventions the Floer chain complex for $P^{-}$ is generated by two gauge equivalence classes isolated flat connections in dimensions 1 and 5 that we shall denote as $\rho_1$ and $\rho_5$ respectively, and we have
\begin{displaymath}
I_*(P^{-}) =\mathbb{Q}_{(1)}\oplus\mathbb{Q}_{(5)}
\end{displaymath}
where the generators over $\mathbb{Z}$ are $\rho_1$ and $\rho_5$.
As shown in \cite{Fy}, $\delta_0(\rho_1)=\pm1$  
and $\delta'_0=0$. On the other hand after an orientation reversal, we have
\begin{displaymath}
I_*(P^{+}) =\mathbb{Q}_{(0)}\oplus\mathbb{Q}_{(4)}
\end{displaymath}
generated by $\rho_0 (=\rho_5)$, $\rho_4 (=\rho_1)$ and $\delta_0=0$, $\delta'_0=\pm\rho_4$. Since $\delta=0$ by the chain rule (\ref{u-chain}) the $u$-map is defined on $I_*(P^{+})$ and \cite{Fy} shows 
\begin{displaymath}
u(\rho_4)=\pm 2\rho_0.
\end{displaymath}

\subsection{Connected sums (i)}\label{subsection-connect}
Here we cite results on connected sums following Fukaya and the simplification due to Donaldson \cite{D}. A further simpification is in \cite{Sca} and we follow that presentation below. Let $Y$, $Y'$ each be admissible, i.e. either a 
$\mathbb{Z}$-homology spheres or such that $w=w_2$ is non-trivial and odd when evaluated on some closed embedded surface. Let $(C_*,\partial)$, $(C'_*,\partial')$ be the respective Floer chain complexes (over $\mathbb{Q}$). We define $I^w_*(Y\sqcup Y')$ by the chain groups  $(\widetilde{C}_*,\widetilde{\partial})$ which are described below. Consider the four complexes
\begin{eqnarray*}
&((C\otimes C')_*, \partial_1) &\partial_1=\partial\otimes I +\epsilon I\otimes\partial', \epsilon =(-1)^{\text{degree}}\\
&(C_*\otimes\mathbb{Q}\langle\theta'\rangle, \partial_2)&\partial_2 =\partial\otimes I\\
&(\mathbb{Q}\langle\theta\rangle\otimes C'_*, \partial_3) &\partial_3=\epsilon\otimes\partial'\\
&((C\otimes C')_*)[3],\partial_4) &\partial_4=-\partial\otimes I -I\otimes\partial'
\end{eqnarray*}
where $\epsilon=(-1)^kI$ on $C_k$ and $(C\otimes C'_*)[3]$ is the complex $(C\otimes C')_*$ but with degrees shifted upwards by 3. If say $Y'$ is admissible but not a homology sphere then the trivial connection $\theta'$ is absent as a generator we would omit the complex $(C_*\otimes\mathbb{Q}\langle\theta'\rangle, \partial_2)$ and a similar statement holds for $Y$.
  $(\widetilde{C}_*,\widetilde{\partial})$ is the direct sum of the four complexes above and with cross terms: 
\begin{displaymath}
\widetilde{\partial}=\partial_1 +\Phi_{1,2}+\Phi_{1,3}+\Phi_{1,4}+ \partial_2+\Phi_{2,4}+\partial_3+\Phi_{3,4}+\partial_4
\end{displaymath}
where
\begin{eqnarray*}
&&\Phi_{1,2}=-I\otimes\widehat{\delta}_{Y'},\quad \Phi_{1,3}= \widehat{\delta}_Y\otimes I,\quad \Phi_{1,4}=2u\otimes I-I\otimes 2u\\
&&\Phi_{2,4}=I\otimes\widehat{\delta}'_{Y'},\quad \Phi_{3,4}=\widehat{\delta}'_Y\otimes I.
\end{eqnarray*}
Here the $\widehat{\delta}$, $\widehat{\delta}'$ maps are obtained from $\delta$, $\delta'$ respectively by identifying $\mathbb{Q}$ with $\mathbb{Q}\langle\theta\rangle$ via $\theta\leftrightarrow 1$. Furthermore in the formulas $\widehat{\delta}$, $\widehat{\delta}'$ should be interpreted as chain maps appropriately. That is, $\widehat{\delta}: C_*\to\mathbb{Q}\langle\theta\rangle$ and $\widehat{\delta}': \mathbb{Q}\langle\theta\rangle\to C_*$ where in the former $\mathbb{Q}\langle\theta\rangle$ is supported in dimension 1 and in the latter dimension 4. 
Let 
\begin{equation}\label{X-cob}
X:Y\sqcup Y'\to Y\# Y'
\end{equation}
denote the cobordism given by the boundary connected sum of the product cobordisms $[0,1]\times Y$ and $[0,1]\times Y'$ where the boundary connected sum is taken near the outgoing ends. Let $C^{\#}_*$ denote the 
Floer chain complex for $I^w_*(Y\# Y')$. $X$ defines a chain map 
\begin{displaymath}
\Phi_{X*}=X_*+ X_{1*}+ X_{2*}+ X^h_*:\widetilde{C}_*\to C^{\#}_*
\end{displaymath}
as follows. The various components are defined by the corresponding 0-dimensional moduli spaces:
\begin{eqnarray*}
&&\langle X_*(\alpha\otimes\beta),\rho\rangle = \#M^{(0)}(\widehat{X})(\alpha, \beta;\rho)\\
&&\langle X_{1*}(\theta\otimes\beta),\rho\rangle = \#M^{(0)}(\widehat{X})(\theta, \beta;\rho)\\
&&\langle X_{2*}(\alpha\otimes \theta),\rho\rangle = \#M^{(0)}(\widehat{X})(\alpha, \theta;\rho)\\
&&\langle X^h_*(\alpha\otimes\beta),\rho\rangle = \#M^{(3)}(\widehat{X})(\alpha, \beta;\rho)\cap\mu(\widehat{\gamma}).
\end{eqnarray*}
The notation $M^{(0)}(\widehat{X})(\alpha, \beta;\rho)$ denotes the $L^2$-moduli space with limits $\alpha, \beta$ at the ingoing ends and $\rho$ and the outgoing end, with the other cases understood similarly. We make some comments about the class $\mu(\widehat{\gamma})$. Here $\gamma$ is a path that connects the boundary components $Y$, $Y'$. $\widehat{\gamma}$ is $\gamma$ (after possibly a reparameterized) extended over the cylindrical-end manifold $\widehat{X}$ as $(t,\text{pt})$, $t\le -1$ in  $(-\infty,0]\times Y$ and as $(-t,\text{pt})$, $t\ge 1$ in $(-\infty,0]\times Y'$. $\mu(\widehat{\gamma})$ is defined to be 1/2 of the degree of the holonomy along $\widehat{\gamma}$ as defined in \cite[\S 7.3.2, \S 7.3.3]{D}. (The factor of $1/2$ is introduced as we use a slightly different normalizations.) In brief, there is a modified holonomy map
\[
H_{\gamma}:M^{(k)}_{\overline{w}}(\widehat{X})(a,b)\to SO(3).
\]
As a divisor, $\mu(\widehat{\gamma})$ is represented by $H_{\gamma}^{-1}(\text{pt})$ for a generic point. 
Note that $\widehat{\gamma}$, being a parameterized curve, is moving in opposite directions along the cylindrical-ends for $Y$ and $Y'$. 
We remark that we are not limited to cobordisms of the form $X$ above. The discussion extends without change to cobordisms
\[
W: Y\sqcup Y'\to Y''
\]
provided the $w_2$ classes involved (if non-trivial) extend over the cobordism, and we can allow $Y''$ to be empty. Then $W$ induces a map
\[
 W_*: I^{w+w'}_*(Y\sqcup Y')\to I^{w''}_*(Y'')
 \]
 where we interpret $I^w(Y'')$ as $\mathbb{Q}$ if $Y''$ is empty. 
 
The following  is essentially Fukaya's result \cite{Fu} as formulated in \cite{Sca}. This formulation covers the various cases where some or all of the $w_2$ classes may be trivial, but in each case the 3-manifold and $w_2$ class is admissible.

\begin{thrm}\label{connectsum}
The chain map $\Phi_{X*}$ induces an isomorphism $\Phi_{X*}: I^{w+w'}_*(Y\sqcup Y')\to I^{w+w'}_*(Y\# Y')$.
\end{thrm}

\subsection{Connected sums (ii)}

Let $P^{-}$ be the Brieskorn homology sphere $B(2,3,5)$ oriented as the boundary of the Milnor fibre, which has intersection form $-E_8$ i.e. we may write $-E_8:\emptyset\to P^-$. Reversing the ingoing and outgoing component we get $-E^*_8: P^+\to\emptyset$ where $P^+$ now has the reversed orientation. We shall follow this orientation notation throughout the paper.

\begin{lem}\label{connect-one}
Let $Y$ be a 3-manifold with $b_1(Y)\neq 0$, and $(Y, w)$ be admissible. Then over $\mathbb{Q}$ 
\begin{displaymath}
I^w_*(Y\# P^{\pm})\cong I^w_*(Y).
\end{displaymath}  
\end{lem}

\proof
We only need to to prove only the case for $P^{+}$ by reversing the orientation of $Y$ if needed. Let the Floer chain complex for $P^{+}$ be $V_*=\mathbb{Q}_{(0)}\oplus\mathbb{Q}_{(4)}$ and the Floer complex for $Y$ be $(C_*,\partial)$. The complex $\widetilde{C}_*$ is then
\begin{displaymath}
(C\otimes V)_*\oplus (C_*\otimes\mathbb{Q}\langle\theta\rangle)\oplus(C\otimes  V)_*[3]
\end{displaymath}
(The 3rd component of the connect sum complex is absent.) The  differentials are
\begin{eqnarray*}
(\partial_1+\Phi_{1,4})(\alpha\otimes \rho_i) &=& (\partial\alpha\otimes\rho_i, 0,2u\alpha\otimes \rho_i - 2\alpha \otimes u\rho_i)\\
(\partial_2+\Phi_{2,4})(\beta\otimes \theta) &=& (0,\partial\beta\otimes \theta,\pm\beta\otimes\rho_4)\\
\partial_4(\alpha\otimes \rho_i) &=& (0,0,-\partial\alpha\otimes\rho_i).
\end{eqnarray*} 
We have a filtration
\begin{displaymath}
0\subset (C\otimes  V)_*[3]\subset (C_*\otimes\mathbb{Q}\langle\theta\rangle)\oplus(C\otimes  V)_*[3]\subset\widetilde{C}_*.
\end{displaymath} 
The $E_2$-page of the spectral sequence is
\begin{displaymath}
(I(Y)\otimes\mbox{\rm coker}\widehat{\delta}'_0)_*[3]\oplus (I_*(Y)\otimes \mbox{\rm ker}\widehat{\delta}'_0)\oplus (I(Y)\otimes V)_*
\end{displaymath}
which simplifies to
\begin{displaymath}
(I_*(Y)\otimes\mathbb{Q}_{(0)})[3]\oplus\{0\}\oplus (I(Y)\otimes V)_*.
\end{displaymath}
The remaining differential is $\Phi_{1,4}: (I(Y)\otimes V)_* \to (I_*(Y)\otimes\mathbb{Q}_{(0)})[3]$ with
\begin{displaymath}
\Phi_{1,4}(\alpha\otimes\rho_0+\beta\otimes\rho_4) = 2(u\alpha\pm 2\beta)\otimes\rho_0.
\end{displaymath}
The map $\alpha\mapsto \alpha\otimes\rho_0\mp u\alpha/2\otimes\rho_4$ defines an isomorphism of $I_*(Y)$ with the kernel of $\Phi_{1,4}$.\qed

Our next task is to find a cobordism $W:Y\# P^+\to Y$ that induces an isomorphism. To this end we shall utilize the cobordism $-E_8^*:P^{+}\to\emptyset$ described above. Recall the definition of $E_8$ as the lattice in $\mathbb{R}^8=\sum_{i=1}^8 a_i e_i$ generated by elements $e_i+e_j$ and $\frac{1}{2}(e_1+\dots e_8)$. Set $w'$ to be the mod~2 class defined by $e_1+e_2+e_3+e_4$. $(Y,w)$ is an admissible pair and the class $w$ clearly extends over the product $[0,1]\times Y$.

\begin{prop}\label{Y-cobor}
Set $W:Y\# P^+\to Y$ to be the cobordism given by the boundary connected sum of the product cobordism $[0,1]\times Y:Y\to Y$ and $-E_8^*:P^{+}\to\emptyset$ where the boundary connected sum it taken at the ingoing ends.
Let $W^{w'}_*$ be the map $I^w_*(Y\# P^+)\to I^w_*(Y)$ where  $w'=e_1+e_2+e_3+e_4$ is as above. Then $W^{w'}_*$ defines an isomorphism (over $\mathbb{Q}$).
\end{prop}

\proof By considering the composition $Z=W\circ X$ where $X$ is as in (\ref{X-cob}), we have a map $Z^{w'}_*: I^w_*(Y\sqcup P^+)\to I^w_*(Y)$. Since $X$ induces an isomorphism according to the connected sum Theorem~\ref{connectsum} it suffices to prove $Z^{w'}_*$ is onto. Observe that we have a connected sum decomposition of $W\circ X$ as $([0,1]\times Y)\#(-E^*_8)$. Working at the chain level, consider a metric on this connected sum with a very long cylindrical tube on the connected sum region (i.e. isometric to $[-T, T]\times S^3$ with the standard product metric) and the behavior of the chain maps as $T\gg 1$. By abuse of notation we shall for the rest of the proof also denote the map induced at the chain level by $W^{w'}_*$, $Z^{w'}_*$ etc. We have
\begin{displaymath}
Z^{w'}_* = W^{w'}_*\circ X_* + W^{w'}_*\circ X_{1*} + W^{w'}_*\circ X^h_*.
\end{displaymath}
The various components are determined by the moduli space on $\widehat{Z}=\widehat{Z}(T)$, the manifold  obtained by adjoining cylindrical tubes to the boundaries. To be specific
\begin{eqnarray*}
\langle (W^{w'}_*\circ X_*)(\alpha\otimes\beta),\rho\rangle &=&\# M^{(0)}_{w+w'}(\widehat{Z})(\alpha, \beta;\rho)\\
\langle (W^{w'}_*\circ X_*)(\alpha\otimes \theta),\rho\rangle &=&\# M^{(0)}_{w+w'}(\widehat{Z})(\alpha, \theta;\rho)\\
\langle (W^{w'}_*\circ X_*)(\alpha\otimes\beta),\rho\rangle &=&\# M^{(3)}_{w+w'}(\widehat{Z})(\alpha, \beta;\rho)\cap\mu(\widehat{\gamma}),
\end{eqnarray*}
where $\gamma$ connects the boundary components $Y$ and $P^+$.
The rest proceeds as in the standard moduli space analysis of a connected sum (see \cite{FM}, \cite{DK}). For $T\gg 1$ the moduli space $M^{(0)}_{w+w'}(\widehat{Z})(\alpha, \beta;\rho)$ must be empty by a dimension count on account of  $w_2$ being non-trival on the connected sum components and $\alpha$, $\beta$ and $\rho$ are irreducible. A similar reasoning shows $M^{(0)}_{w+w'}(\widehat{Z})(\alpha, \theta;\rho)$ is also empty for $T\gg 1$. For $T$ large we have for the last contribution:
\begin{displaymath}
M^{(3)}_{w+w'}(\widehat{Z})(\alpha, \beta;\rho)\cap\mu(\gamma) \cong M^{(0)}_w(\mathbb{R}\times Y)(\alpha,\rho)\times M^{(0)}_{w'}(\widehat{-E^*_8})(\beta)
\end{displaymath}
where we have used the fact that $\mu(\gamma)$ detects the gluing parameter $SO(3)$. For the terms on the right hand side to be zero dimensional spaces we require $\alpha=\rho$, which makes $M^{(0)}_w(\mathbb{R}\times Y)(\alpha,\rho)$ a single point. Secondly we note that only $\beta=\rho_4$ gives a 0-dimensional moduli space on $\widehat{-E^*_8}$ with our choice of $w_2$. This is a relative Donaldson invariant calculation and can be found in \cite[proof of Prop.~2]{Fy}. 
In fact $\# M^{(0)}_{w'}(\widehat{-E^*_8})(\rho_4)=\pm 1$. In the proof of Lemma~\ref{connect-one} it was determined that 
Therefore for $T\gg 1$ the chain map
\begin{displaymath}
Z^{w'}_*(\alpha\otimes \rho_4) =\pm \alpha, \quad \alpha\otimes \rho_4 \in (C_*\otimes V_*)[3].
\end{displaymath}
As it stands $\alpha\otimes\rho_4$ may not be a cycle. However in the proof of Lemma~\ref{connect-one} it was determined that if $\partial\alpha=0$ then $\alpha\otimes \rho_0\mp u\alpha/2 \otimes\rho_4$ is a cycle in the complex for $I^w_*(Y\sqcup P^+)$. Then we see
\[
Z^{w'}_*(\alpha\otimes \rho_0\mp u\alpha/2 \otimes\rho_4) =\pm u\alpha/2.
\]
Since $u$ is an isomorphism, surjectivity follows. By Lemma~\ref{connect-one} this must be an isomorphism.
\qed

Let $n P^{+}$ denotes the $n$-fold connected sum and $-nE^*_8:nP^+\to\emptyset$ the n-fold boundary connected sum of $-E_8^*$. We may iterated the construction of $W:Y\# P^+\to Y$ above to give 
\begin{displaymath}W_n:Y\# nP^+\to Y.
\end{displaymath} 
Repeated application of the Proposition~\ref{Y-cobor} gives the following. 

\begin{cor}\label{cor-cobor}
The map $W^{w'}_{n*}:I^w_{*}(Y\# nP^+)\to I^w_*(Y)$ where  $w'=e_1+e_2+e_3+e_4$ on each $-E^*_8$ factor of $W_n$ is an isomorphism (over $\mathbb{Q}$).
\end{cor}

We need one more result concerning the map on instanton homology induced by the cobordism $W_n$, this time with trivial $w_2$ over the $-E^*_8$ portion.

\begin{prop}\label{zero-cobor}
The map $W^{w'}_{n*}:I^w_{*}(Y\# nP^+)\to I^w_*(Y)$ where  $w'=0$ on each $-E^*_8$ factor is the zero map (over $\mathbb{Q}$).
\end{prop}

\proof When we inspect the limits in the proof of Prop.~\ref{Y-cobor}, the only non-zero term that survives for $T\gg 1$ is
$M^{(0)}_{w+w'}(\widehat{Z})(\alpha, \theta;\alpha)\cong \{\alpha\}\times\{\theta\}$ so $Z^{w'}_*$ is only non-zero on terms of the form $\alpha\otimes \theta$, i.e.
\begin{displaymath}
Z^{w'}_*(\alpha\otimes \theta) =\pm \alpha, \quad \alpha\otimes \theta\in C_*\otimes\mathbb{Q}\langle\theta\rangle.
\end{displaymath}
However in the complex $(\widetilde{C}_*,\widetilde{\partial})$ that computes the instanton homology of the connected-sum, if $\alpha\otimes \theta$ is a cycle then 
$\widetilde{\partial}(\alpha\otimes\theta)=(\partial_2+\Phi_{2,4})(\alpha\otimes \theta) = \partial\alpha\otimes\theta\pm\alpha\otimes\rho_4=0$ and $\alpha$ must be the zero class. \qed

\subsection{Connected sums (iii)}

It was shown by Fukaya \cite{Fu} that $I_*(P^+\#P^+) =\mathbb{Q}^2_{(0)}\oplus \mathbb{Q}^2_{(4)}$. A straightforward application of the connected sum theorem proves the generalization to multiple connected sums.

\begin{lem}\label{nP+}
 $I_*(n P^+) = \mathbb{Q}^n_{(0)}\oplus \mathbb{Q}^n_{(4)}$.
\end{lem} 

\proof As before Floer complex for $I_*(P^+)$ can be expressed as $V_*=\mathbb{Q}_{(0)}\oplus\mathbb{Q}_{(4)}$. Proceed by induction, note the fact that the $\delta_0$ map for $nP^+$ is zero and $\delta'_0$ map is non-zero. This follows from the additivity property of the $h$-invariant  and the fact that $h(P^+) = -1$ \cite{Fy}. Let $C_*$ denote the Floer complex for $nP^+$. The complex $\widetilde{C}_*$ for the connected sum $nP^+\# P^+$ is then 
\begin{displaymath}
(C\otimes V)_*\oplus (C_*\otimes\mathbb{Q}\langle\theta\rangle)\oplus(\mathbb{Q}\langle\theta\rangle\otimes V_*)
\oplus(C\otimes  V)_*[3]
\end{displaymath}
with the non-trivial differentials  
\begin{eqnarray*}
(\partial_1+\Phi_{1,4})(\alpha\otimes \rho_i) &=& (\partial\alpha\otimes\rho_i,0,0, 2u\alpha\otimes \rho_i \pm \alpha \otimes 2u\rho_i)\\
(\partial_2+\Phi_{2,4})(\beta\otimes \theta) &=& (0,\partial\beta\otimes \theta,0,\pm\beta\otimes\rho_4)\\
(\partial_3+\Phi_{3,4})(\theta\otimes\rho_i) &=& (0,0,\theta\otimes\partial\beta,\widehat{\delta}'(\theta)\otimes \rho_i)\\
\partial_4(\alpha\otimes \rho_i) &=& (0,0,0,-\partial\alpha\otimes\rho_i)
\end{eqnarray*} 
The filtration
\begin{displaymath}
0\subset (C\otimes  V)_*[3]\subset (C_*\otimes\mathbb{Q}\langle\theta\rangle)\oplus(\mathbb{Q}\langle\theta\rangle\otimes V_*)\oplus(C\otimes  V)_*[3]\subset\widetilde{C}_*
\end{displaymath} 
determines an $E_2$-page
\begin{displaymath}
(\mbox{\textup{coker}}\widehat{\delta}'_0\otimes\mathbb{Q}_{(0)})[3]\oplus 
(\mbox{\textup{ker}}(\Phi_{2,4}+\Phi_{3,4}))
\oplus (I\otimes V)_*.
\end{displaymath}
$\mbox{\textup{ker}}(\Phi_{2,4}+\Phi_{3,4})\cong\mathbb{Q}$
and is spanned by the dimension 4 class $\mp\widehat{\delta}_0'(\theta)\otimes\theta+\theta\otimes\rho_4$ (recall we should think of $\widehat{\delta}'$ as a chain map appropriately). 
The differential $\Phi_{1,4}:(I\otimes V)_*\to (\mbox{\textup{coker}}\widehat{\delta}_0'\otimes\mathbb{Q}_{(0)})[3]$ is given by
\begin{displaymath}
\Phi_{1,4}(\alpha\otimes\rho_0+\beta\otimes\rho_4) = 2(u\alpha \pm 2\beta)\otimes\rho_0.
\end{displaymath}
The induction hypothesis assumes $I_*=\mathbb{Q}^n_{(0)}\oplus \mathbb{Q}^n_{(4)}$, and since $\delta'\neq 0$ and $\delta'$ maps into dimension 4, we have
\begin{displaymath}
\mbox{\textup{coker}}\widehat{\delta}_0'\cong(\mathbb{Q}^n_{(0)}\oplus \mathbb{Q}^{n-1}_{(4)}).
\end{displaymath}
So when we restrict to dimension zero:
\begin{displaymath}
\Phi_{1,4}: (\mathbb{Q}^n_{(0)}\otimes\mathbb{Q}_{(0)})\oplus (\mathbb{Q}^n_{(4)}\otimes\mathbb{Q}_{(4)})
\to \mathbb{Q}^{n-1}_{(4)}[3]
\end{displaymath}
 which is onto, giving $(\mbox{\rm ker}\Phi_{1,4})_{(0)}\cong \mathbb{Q}^{n+1}$. On the other hand in dimension 4
 \begin{displaymath}
\Phi_{1,4}: (\mathbb{Q}^n_{(0)}\otimes\mathbb{Q}_{(4)})\oplus (\mathbb{Q}^n_{(4)}\otimes\mathbb{Q}_{(0)})
\to \mathbb{Q}^{n}_{(4)}[3]
\end{displaymath}
is onto with $(\mbox{\textup{ker}}\Phi_{1,4})_{(4)}\cong \mathbb{Q}^{n}$. Together with the dimension 4 class $\mbox{\textup{ker}}(\Phi_{2,4}+\Phi_{3,4})\cong\mathbb{Q}$ above we see that $I_*(nP^+\# P^+)\cong\mathbb{Q}^{n+1}_{(0)}\oplus \mathbb{Q}^{n+1}_{(4)}$.\qed

Let $U:nP^{+}\# P^{+}\to nP^{+}$ be the cobordism obtained by taking the boundary connected sum of the product cobordism $[0,1]\times nP^+$ and $-E^*_8: P^+\to\emptyset$, where the boundary connected sum is taken at the ingoing ends. As before, in the $-E^*_8$ portion of the cobordism, we have the $w_2$ class $w'=e_1+e_2+e_3+e_4$.

\begin{prop}
Let $U^{w'}_*$ be the induced map $I_{*}(nP^+\# P^+)\to I_*(nP^+)$ given by the cobordism $U$ and choice of $w_2=w'$. Then $U^{w'}_*$ is of degree zero and surjective.
\end{prop}

\proof The proof follows the same line of argument as Prop.~\ref{Y-cobor}. To see that $U^{w'}_*$ is of degree zero: the moduli space with trivial limits $M_{w'}(\widehat{U})(\theta;\theta)$ has dimension given by the rule $-2(w')^2-3(1-b_1+b^+)\mbox{ \rm mod } 8 \equiv -2(-4)-3\mbox{ \rm mod } 8 \equiv 5 \mbox{ \rm mod } 8$. Then 
by an application of excision for indices of the ASD operator (for example \cite{D}) and definition of the Floer dimension,  we have mod~8:
\begin{displaymath}
\mbox{\textup{dim}} M_{w'}(\widehat{U})(\alpha; \beta) \equiv \mbox{\textup{dim}}(\alpha)+3 + \mbox{\rm dim}M_{w'}(\widehat{U})(\theta;\theta) -\mbox{\textup{dim}}(\beta) \equiv \mbox{\textup{dim}}(\alpha)-\mbox{\textup{dim}}(\beta).
\end{displaymath}
The result follows.

\section{A 4-dimensional characterization of the $h$-invariant for knots}\label{h-char}

Let us review some properties of the $h$-invariant from \cite{Fy}. Fr{\o}yshov constructs, beside the class $\delta_0$, $\delta'_0$ other classes $\delta_{2i}$, $\delta'_{2i}$ by repeated application of the $u$-map i.e. $\delta u^{2i}$, $u^{2i}\delta'$. The classes for $i>0$ are not diffeomorphism invariants but the following holds.  
If $W: Y_1\to Y_2$ is a simply connected negative definite cobordism then the induced maps $W^*$, $W_*$ have the following natural property with respect to the $\delta_{2i}$ and $\delta'_{2i}$. Let $V^{e}_n$, $V^{' e}_n$ denote the spans of $\delta_{2i}$, $\delta'_{2i}$ for $i=0,\dots n$ respectively. 
Then Fr{\o}yshov shows
\begin{equation}\label{h-prop}
W^*(V^{e}_n) = V^{e}_n, \quad W_*(V^{' e}_n) = V^{' e}_n.
\end{equation}
For $n\gg 1$ the increasing sequence of subspaces will stabilize, let us denote this limit as $V^e$, $V^{'e}$. Then the $h$-invariant can be defined as
\begin{displaymath}
h := \mbox{\rm dim}_{\mathbb{Q}} V^e - \mbox{\rm dim}_{\mathbb{Q}} V^{'e}.
 \end{displaymath}
An important observation by \cite{Fy} is that if $h>0$ then $V^{'e}=\{0\}$ and if $h<0$, $V^{e}=\{0\}$. That is to say {\it at least one of $V'$ or $V$  must always be trivial}.

Property (\ref{h-prop}) shows that 
\begin{displaymath}
h(Y_2)\geq h(Y_1)
\end{displaymath}
for simply connected negative definite cobordisms $W: Y_1\to Y_2$. That is to say $h$ is non-decreasing.

Let $W: S^3\to K_{-1}$ be the cobordism that is obtained by adding a 2-handle to the $[0,1]\times S^3$ along $K$, at the boundary $\{1\}\times S^3$, with framing $-1$.  In the following constructions we shall assume $w_2=0$ so we shall omit it from our notation. $W$ is negative define and induces the map $F=W_*: \mathbb{Q}\to I^w_*(K_{-1})$. We can repeat this process taking connected sums with $nP^+$ to give $W_n: nP^+ \to K_{-1}\# nP^+$ (we let $W_0=W$) and degree zero maps ($n>0$)
\begin{equation}\label{Fn-map}
F_{n}=W_{n*}: I_4(nP^+)\to I_4(K_{-1}\# nP^+).
\end{equation}

We can establish a 4-dimensional characterization of the $h$-invariant for the homology spheres obtained by $-1$ surgery on a knot based from the above mentioned properties for the $h$-invariant and our results on stabilization. Recall from Lemma~\ref{nP+} that
\[
I_*(nP^+)\cong\mathbb{Q}^n_{(0)}\oplus \mathbb{Q}^n_{(4)}.
\] 

\begin{prop}\label{stab}
Let $K$ be a knot in $S^3$ and $F_n$ the map (\ref{Fn-map}). Then
\begin{displaymath}
\mbox{\textup{dim ker}}F_{n}=
\left\{
\begin{array}{ll}
n &\mbox{\textup{for $0\leq n < h(K_{-1})$}}\\
h(K_{-1}) &\mbox{\textup{for $n\geq h(K_{-1})$}}.
\end{array}
\right.
\end{displaymath}
Thus $h(K_{-1})$ is the smallest integer $N\geq 0$ such that $\mbox{\rm dim ker}F_{n+1} = \mbox{\rm dim ker}F_{n}$ for $n\geq N$.
\end{prop}

\proof  $h$ is additive under connected sums so we have $h(nP^+)=-n$ which implies $V^{'e}=I_4(nP^+)\cong\mathbb{Q}^n_{(4)}$. Since $K_{-1}$ bounds a simply connected negative definite manifold, $h(K_{-1})=N\geq 0$. We now restrict to  dimension 4, the argument for  dimension 0 is similar. The map $F_n$ takes the form
\begin{displaymath}
F_n: \mathbb{Q}^n_{(4)}\to V^{'e}\subset I_4(K_{-1}\# nP^+).
\end{displaymath}
Assume $n<N$. By additivity $h(K_{-1}\# nP^+) = h(K_{-1}) + h(nP^+)= N-n>0$. We remarked above that we cannot both have $V^e$ and $V^{'e}$ spaces non-trivial therefore $V^{'e}=\{0\}$ and thus $\mbox{\rm dim ker}F_n =n$.  For $n\geq N$ we have $V^{'e}\cong\mathbb{Q}^{n-N}$ and by (\ref{h-prop}), $F_n$ surjects $\mathbb{Q}^n\to \mathbb{Q}^{n-N}$, so 
$\mbox{\rm dim ker}F_{n} =N$.  \qed

\section{Proof of Theorem~\ref{thrm-three}}\label{complete}

Let $W_K: K_0\to\emptyset$ be the defining cobordism for $\delta_K$. We take the boundary connected sum with the cobordism $-nE^*_8:nP^+\to\emptyset$, $n>0$ to get $W_{K,n}:K_0\#nP^+\to \emptyset$. Denote by $\delta_{K,n}$ the homomorphism $I^w_*(K_0\# nP^+)\to \mathbb{Q}$ defined by the zero-dimensional moduli space on $W_{K,n}$, with choice of $w_2=w$ that is equal to $e_1+e_2+e_3+e_4$ in each $-E_8$ factor.

\begin{claim}\label{claim-one-complete}
Under the isomorphism $W^w_{n*}:I^w_*(K_0\# nP^+)\to I^w_*(K_0)$  given in Corollary~\ref{cor-cobor} we have $W^{w*}_n(\delta_K)=\delta_{K,n}$.
\end{claim}
To see the claim observe that we have the factorization
\[
W_{K,n} = W_K \circ W_n.
\]
The cobordism $W_{K,n}$ can be factored another way:
\begin{equation}\label{factor}
W_{K,n} =  (-nE^*_8) \circ U_n, 
\end{equation}
where $U_n:K_0\# nP^+\to nP^+$. We may describe $U_n$ in the following way. Take the defining cobordism $W_K$ and then form the boundary connected sum with the product $[0,1]\times nP^+$. Let $\delta_{-nE^*_8}: I_*(nP^+)\to\mathbb{Q}$ be the zero dimensional Donaldson invariant with the choice of $w_2$ that is equal to $e_1+e_2+e_3+e_4$ in each $-E_8$ factor. Then by the factorization (\ref{factor}) we have
\begin{equation}\label{pull}
U^*_n(\delta_{-nE^*_8}) = \delta_{K,n}.
\end{equation}
Note that $\delta_{-nE^*_8}$ is supported only in dimension 4. 

$U_{n*}$ fits into Floer's exact sequence associated to surgery on $K$, thinking of $K$ as a knot in a fixed 4-ball in $nP^+$:
\begin{equation}\label{exact-seq-fl}
\longrightarrow
{I}_{5}^w(K_{-1}\# nP^+)
\stackrel{}{\longrightarrow}
{I}_{5}^{w}(K_0\# nP^+)
\stackrel{U_{n*}}{\longrightarrow}
{I}_{4}(nP^+)
\stackrel{F_n}{\longrightarrow}
{I}_4^w(K_{-1}\# nP^+)
\longrightarrow.
\end{equation}
Let us now prove the special case $h(K_{-1})=0$ if and only if ${\varphi}^+_e(K)=0$. If $h(K_{-1})=0$ then
 by Prop.~\ref{stab} for $n > 0$, $F_n$ is injective and  by
the exact sequence above $U_{n*}=0$. From the above Claim and (\ref{pull}) $\delta_{K,n}=0$. On the other hand suppose $\delta_{K,n}=0$. From Prop.~\ref{stab}, if $h(K_{-1})>0$ then
$\mbox{\rm dim ker}F_1=1$. In this case it would imply that $U_{1*}$ maps onto $I_4(P^+)\cong\mathbb{Q}$. Consequently $U^*_1(\delta_{-E_8^*}) = \delta_{K,1}\neq 0$ and this would be a contradiction.

\begin{claim}
$U_{n*}$ commutes with the $u$-map action.
\end{claim}

To see the claim: the proof of Theorem~6 in \cite{Fy} determines the behavior of the $u$ map when the cobordism has ends integral homology spheres and is negative definite simply connected. The result applies equally as well to the cobordism $U_n$ since only one end is a integral homology sphere, to give:
\begin{displaymath}
uU_{n*} - U_{n*}u +
\frac{1}{2}\delta'_{nP^+}\circ\delta_{U_n} = 0,
\end{displaymath}
where the term $\delta_{U_n}:I^w_*(K_0\# nP^+)\to\mathbb{Q}$ is defined by counting the zero dimensional moduli space on $\widehat{U}_n$ with trivial connection  limit at the $nP^+$ end. We shall show this term is zero and the claim will follow. In the cobordism $W_{K,n}$, $n>0$ we may make a different choice of $w_2$ on the $-E^*_8$ pieces: choose $w_2=0$ and we have a corresponding invariant $\delta^{\dagger}_{K,n}:I^w_*(K_0\# nP^+)\to\mathbb{Q}$. In the factorization (\ref{factor}) any connection on $\widehat{U}_n$ with trivial limit over the $nP^+$ end may be glued to the trivial connection over $-\widehat{nE^*_8}$ uniquely up to gauge equivalence. Therefore $\delta^{\dagger}_{K,n} = \delta_{U_n}$. Prop.~\ref{zero-cobor} now proves  $\delta_{U_n}=0$.

\begin{claim}
The $n$ vectors $\delta_{-nE^*_8},u^2 \delta_{-nE^*_8}, \dots, u^{2(n-1)}\delta_{-nE^*_8}$ are linearly independent in $I_*(nP^+)^*$.
\end{claim}

This claim follows from applications of Prop.~1 and 2 of \cite{Fy}. Let us recall these results. Suppose we have an integral homology sphere, in our case $nP^+$ bounding a non-standard negative definite simply connected manifold $-nE^*_8$; then by analysis of lowest charge moduli spaces associated the choice of $w_2$, it is possible to determine a relation between Donaldson invariants for  $-nE^*_8$ and the $u$-map. To be specific: a vector $v_0$ in the $-nE_8$ lattice is {\it extremal} if $| v_0^2| \leq |(v_0+2v)^2|$ for all vectors $v$. Given a $w_2$ class $w$ that is extremal, we can consider all the vectors $v$ such that $v^2=w^2$ and $v\equiv w\mbox{\textup{ mod }}2$. Given $w$, \cite[Prop.~1]{Fy} shows how to take a signed count of these vectors $v$  to give a sum denoted $\eta(-nE_8;w)$. Applied to our situation it states
\begin{equation}\label{extremal}
\delta_{-nE^*_8}(u^{i} \delta'_{nP^+}(1)) = 
\left\{
\begin{array}{ll}
0 &\mbox{\rm for $0\leq i<k$}\\
C\eta(-nE_8;w) &\mbox{\rm for $i=k$, $C\neq 0$}
\end{array}
\right.
\end{equation}
where $k=-w^2/2 -1$. In  \cite[Prop.~2]{Fy} it is shown that for  $n=1$ and choice of extremal vector $w_0=e_1+e_2+e_3+e_4$ we have $\eta(-E_8;w_0)\neq 0$ and $k=-w_0^2/2-1=1$. (Note that $\delta'_{P^+}(1)=\pm\rho_4$). For the case $n>1$ choose our $w$ to be the vector $(w_0, w_0, \dots, w_0)$ in $-nE_8$.This remains extremal, since we have a product of definite lattices. This gives $k=2n-1=2(n-1)+1$. 
Futhermore $\eta(-nE_8;w)=\eta(-E_8;w_0)^n\neq 0$. 
Since 
$\delta_{-nE^*_8}(u^i \delta'_{nP^+})=u^i \delta_{-nE^*_8}(\delta'_{nP^+})$ it follows from (\ref{extremal}) that $\alpha=u\delta'_{nP^+}(1)$ satisfies $u^{2j}\delta_{-nE^*_8}(\alpha) = 0$ for $0\le j < n-1$ and $u^{2j}\delta_{-nE^*_8}(\alpha)\neq 0$ for $j=n-1$.
Thus the claimed vectors are independent. 
(As an aside this argument demonstrates that the set of vectors $\{u^{2j}\delta'_{nP^+}(1), 0\le j\le n-1\}$ is independent, which is a demonstration that $h(nP^-)\ge n$.)

To complete the proof:  
Let $n=h(K_{-1})>0$. By the special case we must have ${\varphi}^+_e(K)>0$ and vise-versa.
 By Prop.~\ref{stab} the map $F_n=0$  for $n=h(K_{-1})$ and so by (\ref{exact-seq-fl})  $U_{n*}$ is surjective (and non-trivial). The adjoint $U^*_{n}$ is \textit{injective} and non-trivial and thus 
$U^*_{n}(u^{2i}\delta_{-nE^*_8}), i=0,\dots, n-1$ are also independent in $I^w_*(K_0\# nP^+)^*\cong I^w_*(K_0)^*$ by Claim~\ref{claim-one-complete}. Since $U^*_{n}(u^{2i}\delta_{-nE^*_8})=u^{2i}\delta_{K,n}$ by the same Claim it follows that ${\varphi}^+_e(K)\geq h(K_{-1})$ (assuming $h(K_{-1})>0)$.
On the other hand assume we have $m={\varphi}^+_e(K)>0$ independent vectors $u^{2i}\delta_K, i=0,\dots, m-1$ in $I^w_*(K_0)^*$ and let $n\gg 1$.
Again via Claim~\ref{claim-one-complete} we have the corresponding set of independent vectors $u^{2i}\delta_{K,n}, i=0,\dots, m-1$ in $I^w_*(K_0\# nP^+)^*$.  Since $U^*_{n}(u^{2i}\delta_{-nE^*_8})=u^{2i}\delta_{K,n}$, it follows that $U^*_n$ must be of rank $m$ and thus $U_{n*}$ is also rank $m$. It follows from (\ref{exact-seq-fl}) and Prop.~\ref{stab} that $h(K_{-1})=\mbox{\textup{dim ker}}F_{n}\ge m$. This completes the proof.

\section{Proof of Theorem~\ref{thrm-two}}

Take the connected sum of $Z_K$ with $\overline{\mathbb{CP}^2}$. This does not change the Donaldson invariant as long as we keep $w_2$ zero on the generator of the $H_2(\overline{\mathbb{CP}^2}, \mathbb{Z})$ summand. Take the connected sum of the core of the 2-handle defining $Z_K$ with $\overline{\mathbb{CP}^1}\subset\overline{\mathbb{CP}^2}$ to give surface $C$, with $C\cdot C=-1$. By taking a regular neighbourhood of $B^4\cup C$ we obtain a factorization
\[
Z_K\# \overline{\mathbb{CP}^2} = H\circ Z_{K,-1},
\]
therefore $\psi_K = H_*(\psi_{K,-1})$. We claim that $H_*$ is an isomorphism, and the result will follow. In fact $H_*$ coincides with the isomorphism $b$ in Floer's exact sequence:
\[
\longrightarrow I_*(S^3)=\{0\} \stackrel{a}{\longrightarrow} I_*(K_{-1}) \stackrel{b}{\longrightarrow} I^w_*(K_0)\stackrel{c}{\longrightarrow} I_{*-1}(S^3)=\{0\}\longrightarrow.
\]
To see this, recall that the surgery cobordism $B$ that induces $b$ is defined as follows. Assume $K$ is oriented, and let $(m,l)$ be an oriented meridian-longitude pair on the boundary of $S^3-\nu K^{\circ}$, where $\nu K$ is tubular neighborhood. In $D^2\times S^1$ let $(m',l')$ be defined by $(\partial D^2\times\text{pt}, \text{pt}\times S^1)$. Identify $D^2\times S^1$ with $\nu K$ such that $(m',l')=(m,l)$. Now effect $(-1)$-surgery on $K$ by the new identification $(m',l')=(l-m,-m)$ to give $K_{-1}$.  With this $B$ is the result of attaching a 2-handle to $K_{-1}$ along the core $K'$ of $D^2\times S^1$. We want to regard $K'$ as a knot in $S^3-\nu K^{\circ}$. We accomplish this by isotoping $K'$ off $D^2\times S^1$ and into $S^3-\nu K^{\circ}$. This can be done in different ways, let us make the choice that $K'$ is represented as a $0$-parallel knot to $K$, i.e. the knot that is the result of isotoping $K$ along a Seifert surface for $K$. As such the framing on $K'$ will be the zero-framing as this is the only choice that gives the correct homology for $K_0$, when we think of $K'$ as a knot in the homology sphere $K_{-1}$. This description, after a handle slide, is exactly that  of $Z_K\# \overline{\mathbb{CP}^2}$, where we have attached 2-handles along both $K$ and $K'$ with their given framings. One can verify that $w_2\neq 0$ on $H$. 

\begin{remark}
Let $A$ be the surgery cobordism that induces $a$ above. We have essentially shown that $B\circ A \cong Z_K\# \overline{\mathbb{CP}^2}$. That $b\circ a=0$ but we can have $\psi_K\neq 0$ is explained by the fact that in $B\circ A$ the choice of $w_2$ is such that it is non-zero on the embedded sphere $C$ above but in $Z_K\# \overline{\mathbb{CP}^2}$ the choice is such that $w_2=0$ on $C$.
\end{remark}

\noindent{\textsc{Vallejo, CA}

\end{document}